\documentclass[12pt,reqno]{amsart}

\usepackage[utf8]{inputenc}
\usepackage[T1]{fontenc}
\usepackage{amscd,amssymb,amsmath,amsthm,amsfonts}
\usepackage{calrsfs}
\usepackage{amsaddr}
\usepackage[all]{xy}
\usepackage[british]{babel}
\usepackage[mathscr]{euscript}

\usepackage{stackrel}
\usepackage{mathtools}

\usepackage{tikz}
\usetikzlibrary{matrix,arrows}
\usepackage{tikz-cd}
\usepackage{bm}
\usepackage{orcidlink}
\hypersetup{hypertexnames=false}

\usepackage[papersize={210mm,297mm},lmargin=3cm,rmargin=3cm,top=4cm,bottom=4cm]{geometry}

\theoremstyle{plain}
\newtheorem{theorem}{Theorem}[section]
\newtheorem{corollary}[theorem]{Corollary}
\newtheorem{proposition}[theorem]{Proposition}
\newtheorem{lemma}[theorem]{Lemma}
\theoremstyle{definition}
\newtheorem{definition}[theorem]{Definition}
\newtheorem{example}[theorem]{Example}
\newtheorem{examples}[theorem]{Examples}
\newtheorem{exam-def}[theorem]{Example and Definition}
\theoremstyle{remark}
\newtheorem{remark}[theorem]{Remark}

\def\lam{\lambda}

\def\Ker{\operatorname{Ker}}
\def\Im{\operatorname{Im}}
\def\Coker{\operatorname{Coker}}

\newcommand{\lact}[2]{{}^{#1}\!#2}
\newcommand{\ract}[2]{#1^{#2}}

\newcommand{\K}{\ensuremath{\mathbb{K}}}
\newcommand{\A}{\ensuremath{\mathsf{A}}}
\newcommand{\B}{\ensuremath{\mathsf{B}}}
\newcommand{\id}{\ensuremath{\mathsf{id}}}

\newcommand*{\ldblbrace}{\{\mskip-5mu\{}
\newcommand*{\rdblbrace}{\}\mskip-5mu\}}

\newcommand{\MM}{\ensuremath{\mathsf{M}}}
\newcommand{\II}{\ensuremath{\mathsf{I}}}
\newcommand{\NN}{\ensuremath{\mathsf{N}}}
\newcommand{\ZZ}{\ensuremath{\mathscr{Z}}}
\newcommand{\ZC}{\ensuremath{\mathbf{Z}}}
\newcommand{\T}{\ensuremath{\mathcal{T}}}
\newcommand{\Cat}{\ensuremath{\mathbf{Cat}}}

\newcommand{\Ad}{\ensuremath{\mathsf{Ad}}}

\newcommand{\Ann}{\ensuremath{\mathsf{Ann}}}
\newcommand{\LL}{\ensuremath{\mathcal{L}}}

\newcommand{\RR}{\ensuremath{\mathcal{R}}}

\DeclareMathOperator{\LM}{\mathcal{L}\mathsf{Mul}}
\DeclareMathOperator{\RM}{\mathcal{R}\mathsf{Mul}}
\DeclareMathOperator{\Mul}{\mathcal{B}\ensuremath{\mathsf{Mul}}}
\DeclareMathOperator{\inc}{\ensuremath{\mathsf{inc}}}

\makeatletter
\DeclareFontEncoding{LS2}{}{\noaccents@} \DeclareFontSubstitution{LS2}{stix}{m}{n}
\DeclareSymbolFont{stix@largesymbols}{LS2}{stixex}{m}{n}
\SetSymbolFont{stix@largesymbols}{bold}{LS2}{stixex}{b}{n} 
\DeclareMathDelimiter{\lBrace}{\mathopen} {stix@largesymbols}{"E8} {stix@largesymbols}{"0E}
\DeclareMathDelimiter{\rBrace}{\mathclose} {stix@largesymbols}{"E9} {stix@largesymbols}{"0F}
\makeatother

\numberwithin{equation}{section}

\begin{document}

\title[The center of an associative crossed module]{The center of an associative crossed module}

\author[E. Khmaladze]{Emzar Khmaladze\orcidlink{0000-0001-9492-982X}}
\address{\footnotesize The University of Georgia, Kostava St. 77a, Tbilisi 0171
	Georgia \& A. Razmadze Mathematical Institute of Tbilisi State University
	2, M. Aleksidze II Lane, Tbilisi 0193, Georgia}
\email{e.khmaladze@ug.edu.ge}

\author[M. Ladra]{Manuel Ladra\orcidlink{0000-0002-0543-4508}}
\address{\footnotesize Department of Mathematics \& CITMAga, Universidade de Santiago de Compostela, 15782 Santiago de Compostela, Spain}
\email{manuel.ladra@usc.es}

\author[D. Quiroga]{David Quiroga\orcidlink{0009-0002-8238-752X}}
\address{\footnotesize Department of Mathematics, Universidade de Santiago de Compostela, 15782 Santiago de Compostela, Spain}
\email{david.quiroga@usc.es}

\author[B. Ramos P\'erez]{Brais Ramos P\'erez\orcidlink{0009-0006-3912-4483}}
\address{\footnotesize Department of Mathematics \& CITMAga, Universidade de Santiago de Compostela, 15782 Santiago de Compostela, Spain}
\email{ramosbrais@uniovi.es}


\thanks{The authors were supported by the Agencia Estatal de Investigaci\'on (Spain)  grant PID2024-155502NB-I00 (financial support from MCIN/AEI/10.13039/501100011033 and from FEDER,
	UE)  and by the Xunta de Galicia through the Competitive Reference Groups (GRC), ED431C
	2023/31. Emzar Khmaladze was financially supported  by Shota Rustaveli National Science
	Foundation of Georgia, grant FR-23-271.}
	
\begin{abstract}
In this paper, we define a Drinfeld-type center for an associative crossed module and prove that it is canonically the internal category of an explicit crossed module constructed from Mac Lane’s bimultiplications. The resulting center carries a canonical multiplicative lifting, which induces a natural strict braiding, establishing its correspondence to a braided crossed module. We also examine the relationship between this center and Huq’s categorical center. Finally, we obtain a canonical comparison morphism with Pirashvili’s center of Lie crossed modules via the (braided) Liezation functors, which preserves braidings under an injectivity condition.
\end{abstract}

\maketitle


\section{Introduction}

Whitehead \cite{Wh} introduced the concept of a crossed module of groups as a tool to study relative homotopy groups. Crossed modules of groups are equivalent to strict 2-groups, that is,
 to group objects in the category of small categories.

Later, researchers extensively studied crossed modules and generalized them beyond groups to a range of algebraic categories.
 Crossed modules have been widely employed in homotopy theory, homological algebra, category theory, non-abelian cohomology, and numerous related fields (see, e.g.~\cite{Wa}). 

Gerstenhaber \cite{Ger} developed crossed modules in associative algebras in the context of interpreting the third cohomology group,
 and Dedecker and Lue \cite{DL} did so in the construction of a non-abelian second cohomology theory.

The Lie algebra version was elaborated later by Kassel and Loday \cite{KL}, establishing a framework analogous to the group-theoretic case and highlighting connections to Lie algebra cohomology.  

Within the category of groups, Norrie \cite{Nor} introduced  the center of a crossed module, which coincides with the categorical center in the sense of Huq \cite{Huq}.
Building on this research, Casas and Ladra \cite{CaLa} formulated the theory of the center of a crossed module in the category of Lie algebras.  

Drinfeld introduced the center of a monoidal category (unpublished), and Majid \cite{Ma} and Joyal and Street \cite{JS} developed it. This notion of center categorifies the center of a monoid or a ring.

M.~Pirashvili introduced a Drinfeld-type center for crossed modules of groups \cite{Pim1} and, for crossed modules of Lie algebras, an analogue of the  Drinfeld center for the bracket bifunctor \cite{Pim2}. 

Crossed modules of associative algebras are equivalent to internal categories
in associative algebras, that is, to strict associative 2-algebras. For a
crossed module $(\MM,\A,\mu)$, the associated internal category
$\Cat(\MM,\A,\mu)$ inherits from the multiplication of $\A$ the structure of a strict semimonoidal category \cite{JK}. In
general, it has no unit, since our algebras need not be unital. We introduce
and study a center in this setting, inspired by the Drinfeld center of a
monoidal category and by M.~Pirashvili's centers of crossed modules of groups
and of Lie algebras \cite{Pim1,Pim2}.

Associative crossed modules form a semi-abelian category, so Huq's categorical
notion of center is also available. We compute the Huq center explicitly. We
then show that it embeds into our center through a morphism that, in general, is not a weak equivalence. Finally, we compare our center with Pirashvili's center of Lie
crossed modules via the (braided) Liezation functors. We obtain a canonical
comparison morphism, which preserves braidings under an injectivity condition.

The paper is organized as follows. After this introduction, Section~\ref{S:Prel}  presents preliminaries by
reviewing key concepts related to associative actions, bimultiplications,  associative crossed modules and  crossed squares.
 Section~\ref{S:BCM} revisits the notion of a braided associative crossed module and presents a sufficient condition for the existence of a braiding, namely the existence of a multiplicative lifting.

In Section~\ref{S:center}, we construct the  center of an associative  crossed module $(\MM,\A,\mu)$, denoted by $\ZC(\MM,\A,\mu)$.
 We show that it carries two canonical multiplicative liftings, which induce its braidings (Proposition~\ref{P:liftings}), even when $(\MM,\A,\mu)$ itself admits none.
  We establish several further properties of the center and show that it naturally fits into a homotopy fiber sequence (Corollary~\ref{C:fiber}).
   We also prove that it gives rise to a crossed square associated with $(\MM,\A,\mu)$ (Proposition~\ref{P:square}).
    Section~\ref{S:comp} compares the center of associative crossed modules with Pirashvili's center of Lie
crossed modules~\cite{Pim2}. This comparison is conducted via the (braided) Liezation functors~\cite{FFL1}.
 We construct a canonical comparison morphism and show that it preserves braidings if and only if the two multiplicative liftings of the center coincide, which is the case, for instance, when $\partial$ is injective.

Section~\ref{S:Huq} investigates the relationship between the center introduced here and the Huq center in the semi-abelian category of associative crossed modules.
 Furthermore, it shows that the Huq center embeds into the proposed center via a morphism that, in general, is not a weak equivalence.  

Section~\ref{S:Drin}  establishes that the object constructed in Section~\ref{S:center} corresponds to  the analogue of the  Drinfeld-type center for the multiplication bifunctor
 of the  internal category determined by the crossed module. This result shows that the center is a
braided semimonoidal category corresponding to a braided crossed module. 
In the final section, the center is characterized as the crossed module of elements of $\A$ endowed with coherent null-homotopies.

Throughout the paper we fix a ground field $\K$. All vector spaces and algebras are considered over $\K$, and all linear maps are $\K$-linear maps.


\section{Preliminaries on crossed modules of associative algebras} \label{S:Prel}

An \emph{associative algebra} $\A$ is a $\K$-vector space equipped with a bilinear operation $ \cdot \colon \A\times\A\to\A$, $(a,b) \mapsto a \cdot b =ab$ satisfying the \emph{associative identity}
\begin{equation}\label{Lb_identity}
(ab) c =a(bc),
\end{equation}
for all $a,b,c\in\A$. 
A \emph{morphism of associative algebras} is a linear map which preserves the product, and will be called an \emph{associative morphism}.

A vector subspace $\II$ of an associative algebra $\A$ is called an \emph{ideal} if $ab,\; ba\in \II$ for all $b\in \II$ and $a\in \A$. For example, given an associative algebra $\A$, the subset 
\[
\ZZ(\A)=\Ann(\A) \coloneqq \{a\in\A \mid ab=0=ba  \ \text{for all} \ b\in\A \} \subseteq \A
\]
is an ideal of $\A$ and it is called the \emph{center} (or \emph{annihilator}) of $\A$.

An \emph{abelian associative algebra} is an associative algebra $\A$ with the trivial product, i.e. $ab=0$ for all $a,b\in\A$. It is clear that $\A$ is abelian if and only if $\ZZ(\A)=\A$.

\begin{definition}
Let $\A$ and $\MM$ be associative algebras. An \emph{associative action} of $\A$ on $\MM$ is a couple of bilinear
maps $\A\times \MM\to \MM$, $(a,m)\mapsto {}^am$ and $\MM\times \A\to \MM$, $(m,a)\mapsto m^a$ satisfying the following conditions:
\begin{align}
	 {}^{ab}m &={}^a({}^{b}m), \tag{A1}\label{A1}\\
	m^{ab} & =(m^a)^{b}, \tag{A2}\label{A2}\\
	{}^a (mn) & =({}^a m) n,\tag{A3}\label{A3}\\
	(mn)^a & =m  (n^a), \tag{A4}\label{A4} \\
	({}^a m) ^b & ={}^a(m^{b}), \tag{A5}\label{A5}\\
	(m^a) n & = m ({}^an),\tag{A6}\label{A6}
\end{align}
for all $a,b\in \A$ and $m,n\in \MM$.
\end{definition}
This notion goes back to the extension theory of rings and algebras \cite{Ev,Hoc,Mac} (see also \cite{CDL,LL}), where it appears in terms of bimultiplications (see Definition~\ref{D:bimult} below):
 by \eqref{A3}, \eqref{A4} and \eqref{A6}, for each $a\in\A$ the pair of maps $m\mapsto{}^am$ and $m\mapsto m^a$ is a bimultiplication from $\MM$ to $\MM$, with $\MM$ acting on itself by its product;
  by \eqref{A5}, these bimultiplications are permutable in the sense of \cite{Mac}; and by \eqref{A1} and \eqref{A2}, they depend multiplicatively on $a$. Equivalently,
   a couple of bilinear maps as above is an associative action if and only if the \emph{semidirect product} $\MM\rtimes\A$, that is, the vector space $\MM\oplus\A$ with multiplication
\[
(m,a)(m',a')=(mm'+{}^{a}m'+m^{a'},\,aa'),
\]
is an associative algebra.

For example, if $\MM$ is an  ideal of $\A$, then the product in $\A$ yields an action of $\A$ on $\MM$, that is, ${}^am \coloneqq am$ and $m^a \coloneqq ma$.

Note that the associative action of $\A$ on an abelian associative algebra $\MM$ means that  $\MM$ is just an 
$\A$-bimodule because the axioms \eqref{A3}, \eqref{A4} and \eqref{A6} directly hold.

Consider that $\A$ and $\MM$ are associative algebras with $\A$ acting on $\MM$. We will define the commutator by $ [a,m] \coloneqq {}^a m -m ^a$.
 In the particular case where $\MM=\A$ and the action is exactly the product in $\A$, the commutator is the usual $[a, b]=a b -b a$.

\begin{definition}[\cite{Mac}]\label{D:bimult}
Let $\A$ and $\MM$ be associative algebras with $\A$ acting on $\MM$.

\noindent A \emph{left multiplication} $\lam\colon \A \to \MM$ is a linear map such that 
\begin{equation}\label{Mlemu}
	\lam(ab) =   \lam(a)^b, \quad  a, b \in \A.
\end{equation}
A \emph{right multiplication} $\rho \colon \A \to \MM$ is a linear map such that 
\begin{equation}\label{Rrimu}
	\rho(ab) =	{}^a\rho(b), \quad  a, b \in \A.
\end{equation}
 A \emph{bimultiplication} is a pair $(\lam,\rho) $,  where $\lam\colon \A\to \MM$ is a left multiplication and  $\rho\colon \A\to \MM$ is a right multiplication,  which together satisfy the following relation:
\begin{equation}\label{Mbimul}
	{}^a\lam(b) =    \rho(a)^b, \quad  a, b \in \A.
\end{equation}
\end{definition}

In the particular case where $\MM=\A$ and the action is the product in $\A$, the previous notion of bimultiplication coincides with the usual one introduced by Mac Lane \cite{Mac}
 and with the notion of multiplication described by Hochschild \cite{Hoc}.

We will denote by $\LM(\A, \MM)$, $\RM(\A, \MM)$, and $\Mul(\A, \MM)$ the $\K$-vector spaces of all left multiplications, right multiplications, and bimultiplications  from $\A$ to $\MM$, respectively.

The following lemma follows immediately from the action axioms \eqref{A1}--\eqref{A6}, so we omit its proof. 

\begin{lemma}\label{LmRm_bimul}
Let $\A$ and $\MM$ be associative algebras with $\A$ acting on $\MM$. Let $\LL_m \colon \A \to \MM$ and $\RR_m \colon \A \to \MM$ be the linear maps defined by	
	\begin{equation}\label{maps}
	\LL_m(a) \coloneqq   m^a, \qquad	\RR_m(a) \coloneqq{}^a m,  	 \quad    a \in \A, m \in \MM.
	\end{equation}

The pair  $(\LL_m,\RR_m)$ forms a bimultiplication from $\A$ to $\MM$. Furthermore, the map $\mathcal{B} \colon \MM \to \Mul(\A,\MM)$, $m\mapsto (\LL_m,\RR_m)$, is linear.
\end{lemma}

In particular, for $\MM=\A$ with the action given by the product, the left multiplication map $\LL_a \coloneqq a \cdot \colon  \A \to \A$, $b \mapsto ab$, 
and  the  right multiplication map $\RR_a \coloneqq \cdot a  \colon  \A \to \A$, $b \mapsto ba$, together form a bimultiplication $(\LL_a,\RR_a)$  from $\A$ to $\A$.

In the following definition, we introduce the main object of this article: crossed modules of associative algebras.
\begin{definition}[\cite{DL,Ger}]
	A \emph{crossed module} $(\MM,\A,\mu)$ is an associative morphism $\mu \colon \MM\to\A$ together with an associative action of $\A$ on $\MM$ such that
\begin{align}
\mu({}^am)=a\mu(m), \quad \mu(m^a)=\mu(m)a, \quad &  \  (\mu \  \text{is}  \ \A\text{-equivariant}) \label{CM1} \tag{CM1}\\
^{\mu(m)}n=m n=m^{\mu(n)}    \quad &  \  (\text{Peiffer identity}) \label{CM2}\tag{CM2}
\end{align}
for all $a\in \A$, $m, n\in \MM$.

An associative crossed module is said to be \emph{abelian} if both $\MM$ and $\A$ are abelian associative algebras and the action of $\A$ on $\MM$ is trivial.

Given another associative crossed module $(\NN,\B,\nu)$, a \emph{morphism of associative crossed modules} 
$(f,g) \colon  (\MM,\A,\mu) \to (\NN,\B,\nu)$ is a pair of associative morphisms, $f \colon  \MM \to \NN$  and 
$g \colon  \A \to \B$, such that 
\begin{equation}\label{morasoxmod}
	 \nu f = g \mu, \quad f(^am)=  {}^{g(a)}f(m), \quad f(m^a)=  f(m)^{g(a)},  \quad a \in \A, m \in \MM.
\end{equation}
\end{definition}

Given a crossed module  $(\MM,\A,\mu)$, the image  $\Im(\mu)$ is an  ideal of $\A$  
and the kernel $\Ker(\mu)$ is an $\A$-invariant  ideal contained in the center of $\MM$. Thus, $\Ker(\mu)$ is an abelian associative algebra and it inherits a structure of an $\A$-bimodule;
 since $\Im(\mu)$ acts trivially on $\Ker(\mu)$ by \eqref{CM2}, it is in fact a $\Coker(\mu)$-bimodule.

We denote by $\pi_0(\MM,\A,\mu)\coloneqq\Coker(\mu)$ and $\pi_1(\MM,\A,\mu)\coloneqq \Ker(\mu)$ the homotopy groups of the crossed module $(\MM,\A,\mu)$. Recall that a morphism of crossed modules
$(f,g) \colon  (\MM,\A,\mu) \to (\NN,\B,\nu)$
is a \emph{weak equivalence} whenever the induced maps $\pi_1(\MM,\A,\mu) \to \pi_1(\NN,\B,\nu)$ and $\pi_0(\MM,\A,\mu) \to \pi_0(\NN,\B,\nu)$ are isomorphisms.
 
 \begin{examples} \label{E:CM} \hfill
 	\begin{enumerate}
 		\item Let $\A$ be an associative algebra and $\II$ an  ideal of $\A$.  The inclusion map  $\II \hookrightarrow \A$ defines a crossed module,
 		 where $\A$ acts on $\II$ via the product operation. As a special case, we have the crossed modules $ \id \colon\A \to \A$ and
 		$ 0 \colon 0 \to \A$.
 		\item  If $\MM$ is an $\A$-bimodule, then  the zero map $ 0 \colon \MM \to \A$ forms a crossed module.
 		\item Let $\A$ be the two-dimensional associative algebra with basis $\{e,f\}$ whose only non-zero products of basis elements are
 		\[
 		e^2=e, \qquad ef=f.
 		\]
 		 The algebra $\A$ is non-commutative and is non-unital. Moreover, $\ZZ(\A)=0$, $[\A,\A]=\K f$, and $\K f$ is the only proper non-zero ideal of $\A$.

 		Let $\MM$ be the two-dimensional abelian associative algebra with basis $\{u,v\}$, and define the action of $\A$ on $\MM$ by
 		\[
 		{}^e u=u,
 		\]
 		all the other actions of basis elements being zero. Then the linear map
 		\[
 		\mu \colon \MM\longrightarrow \A, \qquad \mu(u)=f, \quad \mu(v)=0,
 		\]
 		is an associative crossed module. 
 	\end{enumerate}
 \end{examples}
 
 The following result provides an equivalent characterization of crossed modules.
 \begin{proposition}
 	Let $\A$ be an associative algebra.	A linear map  $\mu \colon \MM\to\A$ is an associative crossed module if and only if $\MM$ is an $\A$-bimodule and $\mu$ is an $\A$-equivariant morphism (i.e.  
 	$\mu ({}^a m)=a  \mu(m)$, $\mu (m^a)=\mu(m) a$)  satisfying the following identity:
 	\begin{equation}\label{char_acm}
 	{}^{\mu(m)} n = m^{\mu(n)}
 	\end{equation}
 	for all $a\in\A$, $m, n\in \MM$.	
 \end{proposition}
 \begin{proof}``$\Longrightarrow$'' It is direct by the axioms of an associative crossed module.
 	
 	``$\Longleftarrow$'' Under these hypotheses, $\MM$ is an associative algebra with the product $mn\coloneqq{}^{\mu(m)}n\stackrel{\eqref{char_acm}}{=}m^{\mu(n)}$. Moreover, since $\mu$ is $\A$-equivariant, the linear map $\mu$ is an associative morphism. Then, $(\MM,\A,\mu)$ is an associative crossed module with the action of $\A$ on $\MM$ induced by the $\A$-bimodule structure.
 \end{proof}
 
 Moreover, given an associative crossed module $(\MM,\A,\mu)$, the following result shows that the vector space of bimultiplications of $\A$ in $\MM$ becomes an associative algebra.
\begin{lemma}\label{L:bimul}
 If $(\MM,\A,\mu)$ is an associative crossed module, then $\Mul(\A, \MM)$ is an associative algebra under the product
\begin{equation}\label{E:bimul}
 (\lam,\rho) \cdot (\lam',\rho')  \coloneqq  (\lam \mu \lam',  \rho' \mu \rho).
\end{equation}
\end{lemma} 
\begin{proof}
	It is straightforward to show that the product is bilinear and \eqref{Lb_identity} follows from the associativity of composition. Therefore, it only remains to prove that the product is well-defined. Indeed, $\lambda\mu\lambda'\in \LM(\A,\MM)$ because
	\[
		\lambda\mu\lambda'(ab)\stackrel{\eqref{Mlemu}\text{ for $\lambda'$}}{=}\lambda\mu(\lambda'(a)^{b})\stackrel{\eqref{CM1}}{=}\lambda(\mu(\lambda'(a))b)\stackrel{\eqref{Mlemu}\text{ for $\lambda$}}{=}\lambda\mu\lambda'(a)^b.
	\]
	Analogously, $\rho'\mu\rho\in\RM(\A,\MM)$, and the compatibility \eqref{Mbimul} is obtained as follows:
	\[
		{}^a(\lambda\mu\lambda'(b))\stackrel{\eqref{Mbimul}\text{ for $(\lambda,\rho)$}}{=}\rho(a)^{\mu\lambda'(b)}\stackrel{\eqref{CM2}}{=}\rho(a)\lambda'(b)\stackrel{\eqref{CM2}}{=}{}^{\mu\rho(a)}\lambda'(b)\stackrel{\eqref{Mbimul}\text{ for $(\lambda',\rho')$}}{=}\rho'\mu\rho(a)^{b}.
	\]
\end{proof}
 
To conclude this preliminary section, we recall the notion of a crossed square.
\begin{definition}[\cite{El}]\label{square} A \emph{crossed square} of associative algebras is a commutative diagram of associative algebras
	\[	
	\begin{tikzcd}
		\MM \arrow[d,"\mu_1",swap] \arrow[r, "\mu_2"]  & \NN \arrow[d,"\nu_2"] \\
		\A \arrow[r,"\nu_1",swap]                     & \B           
	\end{tikzcd} 
	\]
together with actions of  $\B$ on $\MM, \A, \NN$ (consequently, $\A$ acts on $\MM$ and on $\NN$ via $\nu_1$, and  $\NN$ acts on $\MM$ and on $\A$ via $\nu_2$), and
 two bilinear maps $h_1 \colon \A \times \NN \to \MM$, $h_2 \colon \NN  \times \A \to \MM$, that satisfy the following properties:
 \begin{equation}\tag{SC1}  \label{SC1}
\begin{aligned}
&	\text{The morphisms} \  \mu_1,\mu_2,\nu_1,\nu_2, \nu_1\mu_1  \ \text{and} \  \nu_2\mu_2 \ \text{are crossed modules};\\
&	\text{The morphisms} \ \mu_1 \  \text{and} \  \mu_2 \ \text{preserve the actions of} \ \B \ \text{on} \ \MM; 
\end{aligned}
\end{equation}
\begin{equation}\tag{SC2} \label{SC2}
	\begin{aligned}
		\mu_1 h_1(a,n)= a^n,& \qquad \qquad  \mu_2 h_1(a,n)=  {} ^a{}n , \\
		\mu_1 h_2(n,a)=   {}^n{}a,&  \qquad \qquad \mu_2 h_2(n,a)= {} n^a;
	\end{aligned}
\end{equation}
\begin{equation}\tag{SC3} \label{SC3}
	\begin{aligned}
		 h_1(\mu_1(m),n)= m^n,&   \qquad \qquad  h_1(a,\mu_2(m))= {}^a{}m,\\
	 h_2(n,\mu_1(m))= {}^n{}m,&  \qquad \qquad h_2(\mu_2(m),a)= m^a;
	\end{aligned}
\end{equation}
\begin{equation}\tag{SC4} \label{SC4}
	\begin{aligned}
	 h_1(a a',n)=   {}^a{}h_1(a',n),&   \qquad  \qquad  h_1(a,nn')= h_1(a,n)^{n'},  \\
  h_2(n,a a')= h_2(n,a)^{a'} ,& \qquad   \qquad 	 h_2(nn',a)=  {}^n{}h_2(n',a),\\
  {}^n{}h_1(a,n')    =  h_2(n,a)^{n'},& \qquad  \qquad	h_1(a,n)^{a'}  = {}^a{}h_2(n,a');\\
\end{aligned}
\end{equation}
\begin{equation}\tag{SC5} \label{SC5}
	\begin{aligned}
h_1(a,n^b)= h_1(a,n)^b,&  \qquad  \qquad  h_1({}^b{}a,n) = {}^b{} h_1(a,n) ,   \\
	h_2({}^b{}n,a)  = {}^b{} h_2(n,a),&  \qquad 	  \qquad    h_2(n,a^b) = h_2(n,a)^b,  \\
 h_1(a,{}^b{}n)= h_1(a^b,n),& \qquad  \qquad  h_2(n^b,a) = h_2(n,{}^b{}a);
\end{aligned}
\end{equation}
for all $ m \in \MM,  \, a,a'\in \A, \,  n,n' \in \NN$,  and $ b \in \B $.
\end{definition}

\begin{example} Let $\B$ be an associative algebra, $\A$, $\NN$  any two ideals of $\B$ and $\MM = \A \cap \NN$. If we define
	\begin{align*}
		h_1 \colon \A \times \NN \to \A \cap \NN, \ h_1(a,n) \coloneqq a \,n,\\
		h_2 \colon \NN  \times \A \to \A \cap \NN, \ h_2(n,a) \coloneqq n \, a,
	\end{align*}
	then, 
	\[	
	\begin{tikzcd}
		\A \cap \NN \arrow[d,"\mu_1", hook, swap] \arrow[r, hook,"\mu_2"]  & \NN \arrow[d,hook,"\nu_2"] \\
		\A \arrow[r,hook,"\nu_1",swap]                     & \B           
	\end{tikzcd} 
	\]
is a crossed square of associative algebras, where $\mu_1,\mu_2,\nu_1,\nu_2$ are the inclusion maps.
\end{example}


\section{On braided  associative crossed modules}\label{S:BCM}

\begin{definition}[{\cite[Definition~3.3]{FFL1}}]
A \emph{braiding} on an associative crossed module  $(\MM,\A,\mu)$ is a  
	bilinear map $\{-,-\} \colon \A \times \A\xrightarrow{}\MM$, $(a,b)\mapsto \{a,b\}$, satisfying:
	\begin{align}
		\mu\{a,b\} & =[a ,b],\tag{BA1}\label{BA1}\\
		\{\mu(m),a \} & =-[a,m],\tag{BA2}\label{BA2}\\
		\{a,\mu (m)\} & = [a,m],\tag{BA3}\label{BA3}\\
		\lbrace a b ,c \rbrace  & =  {}^a \lbrace b,c\rbrace + \lbrace a,c\rbrace ^b ,\tag{BA4}\label{BA4} \\
		\lbrace a  , bc \rbrace  & =  {}^b \lbrace a,c\rbrace + \lbrace a,b\rbrace ^c ,\tag{BA5}\label{BA5} 
	\end{align}
	for all $a,b,c\in \A$, $m\in \MM$.
	In this case, we say that $(\MM,\A,\mu)$ is a \emph{braided associative crossed module}.
\end{definition}

\begin{remark}
	The definition of a braided associative crossed module in \cite[Definition~3.3]{FFL1}  contains one supplementary axiom:
	\[ \{\mu (m),\mu (n) \}=[m,n],\quad m,n\in \MM.\]
	This axiom is redundant because
it follows from \eqref{BA2} and the Peiffer identity  \eqref{CM2}:
	\[\{\mu(m), \mu (n)  \} \stackrel{\text{\eqref{BA2}}}{=} - [\mu (n),m]= -({}^{\mu (n)} m - m^{\mu (n)}) \stackrel{\eqref{CM2}}{=} -(nm -mn) =[ m,n] .\]
\end{remark}

\begin{examples} \label{E:braidings} \hfill
	\begin{enumerate}
	\item In the case of the identity crossed module $(\A,\A,\id)$,
	the unique braiding is given by $\{a,b\} =  [a,b]$.
	
	\item  Let $\A=\K e$ be the one-dimensional abelian associative algebra, i.e. $e^2=0$, and let $\MM=\K v$ be an $\A$-bimodule with trivial action  of $\A$ on $\MM$. 
	In addition to the trivial braiding $\{e,e\}=0$ on $(\MM,\A,0)$, non-trivial braidings also exist: for every $k \in \K$, the bilinear map determined by
	$\{e,e\}=kv$ is a braiding on $(\MM,\A,0)$.

	\item Consider the crossed module $(\MM,\A,\mu)$ of Example~\ref{E:CM}~(3). It admits a unique braiding, namely the bilinear map determined by
	\[
	\{e,f\}=u, \qquad \{f,e\}=-u, \qquad \{e,e\}=\{f,f\}=0.
	\]

	A direct check on basis elements shows that this bilinear map satisfies \eqref{BA1}--\eqref{BA5}. Let us show the uniqueness. Applying \eqref{BA3} and \eqref{BA2} with $m=u$, we obtain $\{e,f\}=\{e,\mu(u)\}={}^eu-u^e=u$, $\{f,f\}=0$ and $\{f,e\}=\{\mu(u),e\}=u^e-{}^eu=-u$. By \eqref{BA1}, $\{e,e\}\in\Ker(\mu)=\K v$, and then \eqref{BA4} gives $\{e,e\}=\{ee,e\}={}^e\{e,e\}+\{e,e\}^e=0$. 

	\item Let $\A$ be the algebra of Example~\ref{E:CM}~(3). The crossed module $(0,\A,0)$ does not admit a braiding, since \eqref{BA1} would give $0=\mu\{e,f\}=[e,f]=f$.
	\end{enumerate} 
\end{examples} 

\begin{definition}\label{D:ml}
	A \emph{multiplicative lifting} on an associative crossed module $(\MM,\A,\mu)$ is a bilinear map
	$\langle-,-\rangle\colon \A\times \A\to \MM$ such that
	\begin{align}
		\mu \langle a,b\rangle & = ab,\label{ML1} \tag{ML1}\\
		\langle a,\mu (m)\rangle & = {}^a m, \label{ML2} \tag{ML2}\\
		\langle \mu (m),a\rangle & = m^a, \label{ML3} \tag{ML3}\\
		\langle ab,c\rangle & =\langle a,bc\rangle, \label{ML4} \tag{ML4}
	\end{align}
	for all $a,b,c\in\A, m\in\MM$.
\end{definition}

Applying \eqref{ML3} with $a=\mu n$ and using \eqref{CM2} yields $\langle\mu (m),\mu (n)\rangle=mn$. Therefore, a multiplicative lifting simultaneously lifts all four multiplications: that of $\A$, that of $\MM$, and the two actions.

A crossed module with a multiplicative lifting can be characterized equivalently as given in the following proposition.

\begin{proposition}\label{P:single}
An associative crossed module equipped with a multiplicative lifting can equivalently be defined as a linear map of vector spaces $\mu \colon \MM\to \A$
 together with a bilinear map $\langle-,-\rangle\colon \A \times \A\to \MM$ satisfying the following single axiom:
\begin{equation} \label{E:pas}
	\bigl\langle \mu \langle a,b\rangle,\,c\bigr\rangle
	=\bigl\langle a,\,\mu \langle b,c \rangle\bigr\rangle ,
\end{equation}
	for all $a,b,c\in \A$.
\end{proposition}

\begin{proof}
	If $\langle-,-\rangle$ is a multiplicative lifting then \eqref{E:pas} is \eqref{ML4}
	rewritten by \eqref{ML1}. 
	
	Conversely, let $(\mu,\langle-,-\rangle)$ be a pair satisfying \eqref{E:pas} and define
	\begin{align*}
	a b & \coloneqq \mu \langle a ,b \rangle,\qquad  
	m n  \coloneqq \langle \mu (m),\mu (n) \rangle, \\
	{}^a m & \coloneqq  \langle a,\mu(m)\rangle,\qquad 
	m^a   \coloneqq  \langle \mu(m),a\rangle,
	\end{align*}
	for all $a,b\in \A$ and $m,n\in\MM$. With the products defined in this way, $\A$ and $\MM$ are associative algebras, and $\mu$ is automatically an associative algebra morphism.
	 Moreover, with the actions defined as above, it is straightforward to verify that $(\MM,\A,\mu)$ is an associative crossed module with multiplicative lifting $\langle -,-\rangle$.
\end{proof}

In the following result we will show that every multiplicative lifting on an associative crossed module induces a braiding.
\begin{proposition}\label{P:ml-braiding}
	If $\langle-,-\rangle$ is a multiplicative lifting on $(\MM,\A,\mu)$, then
	\[
	\{a,b\} \coloneqq \langle a,b \rangle-\langle b,a\rangle
	\]
	is a braiding on $(\MM,\A,\mu)$.
\end{proposition}
\begin{proof} Let us prove the axioms \eqref{BA1}--\eqref{BA5}.
	\begin{itemize}
		\item[\eqref{BA1}:]$\mu\{a,b\}=\mu\langle a,b \rangle-\mu\langle b,a \rangle\stackrel{\eqref{ML1}}{=}ab-ba=[a,b]$.
		\item[\eqref{BA2}:]$\{\mu(m),a\}=\langle\mu(m),a\rangle-\langle a,\mu(m)\rangle\stackrel{\eqref{ML2},\eqref{ML3}}{=}m^{a}-{}^{a}m=-[a,m]$.
		\item[\eqref{BA4}:]
		\begin{align*}
			{}^{a}\{b,c\}+\{a,c\}^{b} & = {}^{a}\langle b,c\rangle -{}^{a}\langle c,b\rangle+\langle a,c\rangle^{b}-\langle c,a\rangle^{b}\\
			&\stackrel{\eqref{ML2},\eqref{ML3}}{=} \langle a,\mu\langle b,c\rangle\rangle-\langle a,\mu\langle c,b\rangle\rangle+\langle\mu\langle a,c\rangle,b\rangle -\langle\mu\langle c,a\rangle,b\rangle \\
			& \stackrel{\eqref{ML1}}{=} \langle a,bc\rangle-\langle a,cb\rangle+\langle ac,b\rangle -\langle ca,b\rangle\\
			& \stackrel{\eqref{ML4}}{=} \langle ab,c\rangle -\langle c,ab\rangle=\{ab,c\}.
		\end{align*} 
	\end{itemize}
	The axioms \eqref{BA3} and \eqref{BA5} are analogous to \eqref{BA2} and \eqref{BA4}, respectively.
\end{proof}

\begin{remark}\label{R:noML}
	The crossed module $(\MM,\A,\mu)$ of Example~\ref{E:CM}~(3) admits no multiplicative lifting, since \eqref{ML1} would give $\mu\langle e,e\rangle=e\notin\Im(\mu)=\K f$.
	 Hence its braiding, described in Example~\ref{E:braidings}~(3), does not arise from Proposition~\ref{P:ml-braiding}.
\end{remark}


\section{The center of an associative crossed module}\label{S:center}

The aim of this section is to construct, from each crossed module, a braided crossed module, which will be referred to as its center.

Let $(\MM,\A,\mu)$ be an  associative crossed module, and consider the set of triples
\[
\T(\MM,\A,\mu) \coloneqq
\{(\lam, \rho,a) \mid a\in \A,\ (\lam, \rho)\in\Mul(\A,\MM) \ \text{satisfying \eqref{Z1}--\eqref{Z4}}\},
\]
where the conditions \eqref{Z1}--\eqref{Z4} are the following identities:
\begin{align}
	\mu\lam(b) & = ab,  \label{Z1} \tag{Z1}\\
	\mu\rho(b) & = ba, \label{Z2} \tag{Z2}\\
	\lam\mu(m) & = {}^a  m,  \label{Z3} \tag{Z3}\\
	\rho\mu(m) & = m^a , \label{Z4} \tag{Z4}
\end{align}
for all $b\in \A$, $m\in \MM$.	

Clearly,  $\T(\MM,\A,\mu)$ is a $\K$-vector space.
\begin{lemma}\label{partial_linear}
	Let $(\MM,\A,\mu)$	be an  associative crossed module. For every $m \in \MM$, the triple $(\LL_m,\RR_m,\mu(m))$ belongs to $\T(\MM,\A,\mu)$.
	 Consequently, $ \partial \colon \MM \longrightarrow  \T(\MM,\A,\mu)$  is a linear map defined by 
		\begin{equation}\label{def_partial}
			\partial (m) \coloneqq (  \LL_m,\RR_m, \mu(m)).
		\end{equation}		
\end{lemma}

\begin{proof} Note that Lemma \ref{LmRm_bimul} asserts that $(\LL_m, \RR_m)\in\Mul(\A,\MM)$. Moreover, the fact that $\mu$ is $\A$-equivariant \eqref{CM1} implies that \eqref{Z1} and \eqref{Z2} hold, whereas \eqref{Z3} and \eqref{Z4} follow by the Peiffer identity \eqref{CM2}.
	\end{proof}
	
The previous lemma suggests that the natural candidate for a braided crossed module constructed from $(\MM,\A,\mu)$ is $(\MM,\T(\MM,\A,\mu),\partial)$. To establish this, we first show that $\T(\MM,\A,\mu)$ can be endowed with an associative algebra structure in such a way that $\partial$ becomes an associative morphism.		
\begin{lemma}\label{L:comp}
Let $(\MM,\A,\mu)$ be an associative crossed module.
	
For all
	$(\lam,\rho,a),(\lam',\rho',b)\in\T(\MM,\A,\mu)$, we have
	\[\partial\lam(b)=(\LL_{\lam(b)},\RR_{\lam(b)}, ab) = (\LL_{\rho'(a)},\RR_{\rho'(a)},ab) = \partial\rho'(a). \]
\end{lemma}
\begin{proof}
	We have to show that $\LL_{\lambda(b)}=\LL_{\rho'(a)}$ and $\RR_{\lambda(b)}=\RR_{\rho'(a)}$. Let $c\in \A$.
	
	\fbox{$\LL_{\lambda(b)}=\LL_{\rho'(a)}$}  We have that:
	\begin{align*}
		\LL_{\lambda(b)}(c) & =\lambda(b)^{c}\stackrel{\eqref{Mlemu}\text{ for $\lambda$}}{=}\lambda(bc)\stackrel{\eqref{Z1}\text{ for $(\lambda',\rho',b)$}}{=}\lambda\mu\lambda'(c)\stackrel{\eqref{Z3}\text{ for $(\lambda,\rho,a)$}}{=}{}^{a}\lambda'(c)\\
		& \stackrel{\eqref{Mbimul}}{=}\rho'(a)^{c}=\LL_{\rho'(a)}(c).
	\end{align*}
	
	\fbox{$\RR_{\lambda(b)}=\RR_{\rho'(a)}$} Following similar arguments as above we have that:
	\begin{align*}
		\RR_{\lambda(b)}(c) & = {}^{c}\lambda(b)\stackrel{\eqref{Mbimul}}{=}\rho(c)^{b}\stackrel{\eqref{Z4}\text{ for $(\lambda',\rho',b)$}}{=}\rho'\mu\rho(c)\stackrel{\eqref{Z2} \text{ for $(\lambda,\rho,a)$}}{=}\rho'(ca)\\
		& \stackrel{\eqref{Rrimu}}{=} {}^{c}\rho'(a)=\RR_{\rho'(a)}(c). \qedhere
	\end{align*}
\end{proof}

\begin{remark}
	Note that Lemma~\ref{L:comp} does not assert that  $\lam(b)=\rho'(a)$. By Lemma~\ref{L:comp}, the difference $\lam(b)-\rho'(a)$ always lies in $\Ker(\partial)$, so $\lam(b)=\rho'(a)$ whenever $\partial$ is injective. In general, however, $\lam(b)$ and $\rho'(a)$ may be different.
	
	 Let $\A=\K e$ be the one-dimensional associative algebra with $e^2=0$, and  $\MM=\K v$  the one-dimensional abelian associative algebra, regarded as an $\A$-bimodule with trivial action. By Example~\ref{E:CM}~(2), $(\MM,\A,0)$ is an associative crossed module. Since $\mu=0$ and the product of $\A$ and the action of $\A$ on $\MM$ are trivial,  $\T(\MM,\A,0)$ consists of all the triples $(\lam,\rho,a)$ with $\lam,\rho\colon\A\to\MM$ linear maps and $a\in\A$.
	
	Let $\lam \colon\A\to\MM$ be the linear map given by $\lam(e)=v$, and take $(\lam,0,0)$ and $(0,0,e)$ in $\T(\MM,\A,0)$, i.e. $\rho'=0$. Then
	\[
	\lam(e)=v\neq0=\rho'(0),
	\]
	whereas $\partial\lam(e)=\partial(v)=(\LL_v,\RR_v,0)=(0,0,0)=\partial\rho'(0)$, since $v$ is annihilated by $\A$.
\end{remark}

\begin{proposition}\label{prop:algebra}
	Let $(\MM,\A,\mu)$	be an  associative crossed module. The vector space $\T(\MM,\A,\mu)$ has the structure of an associative algebra with the product defined by
	\begin{equation}\label{eq:prod}
		 (\lam,\rho,a) \cdot (\lam',\rho',b)  \coloneqq \partial \lam(b)=\partial \rho'(a).
	\end{equation}
	 Explicitly, $(\lam,\rho,a) \cdot (\lam',\rho',b)=(\lam \mu \lam',  \rho' \mu \rho, ab)$; that is, $\T(\MM,\A,\mu)$ is a subalgebra of the product of the algebra $\Mul(\A, \MM)$ of Lemma~\ref{L:bimul} and $\A$.
\end{proposition}
\begin{proof}
	The explicit formula follows from the computations in the proof of Lemma~\ref{L:comp}, which show that $\LL_{\lam(b)}=\lam\mu\lam'$ and $\RR_{\lam(b)}=\RR_{\rho'(a)}=\rho'\mu\rho$. The fact that the product is bilinear is a direct consequence of the linearity of $\partial$ (see Lemma \ref{partial_linear}). Hence, it only remains to show the associativity. On the one hand,
	\begin{align*}
		[(\lam,\rho,a)(\lam',\rho',b)](\lam'',\rho'',c) & = \partial\lam(b)(\lam'',\rho'',c)=(\LL_{\lam(b)},\RR_{\lam(b)}, ab)(\lam'',\rho'',c)\\
		& = \partial\LL_{\lambda(b)}(c)=\partial(\lambda(b)^{c})\stackrel{\eqref{Mlemu}}{=}\partial\lambda(bc),
	\end{align*}
	while, on the other hand,
	\[(\lam,\rho,a)[(\lam',\rho',b)(\lam'',\rho'',c)]=(\lam,\rho,a)\partial\lambda'(c)=(\lam,\rho,a)(\LL_{\lam'(c)},\RR_{\lam'(c)}, bc)=\partial\lambda(bc).\]
	Thus, we conclude that the product is associative.
\end{proof}

\begin{lemma}\label{L:omega_mor}
	Let $(\MM,\A,\mu)$	be an  associative crossed module. The map 
	\[\omega \colon \T(\MM,\A,\mu) \xrightarrow{\quad} \A, \quad   \omega(\lam,\rho,a) \coloneqq a,\]
is an  associative morphism.
\end{lemma}	
\begin{proof}
	The map $\omega$ is clearly linear and it also preserves the product. Indeed,
	\[
		\omega((\lambda,\rho,a)(\lambda',\rho',b))=\omega(\partial\lambda(b))=\omega((\LL_{\lam(b)},\RR_{\lam(b)}, ab))=ab=\omega(\lambda,\rho,a)\omega(\lambda',\rho',b).
	\]
\end{proof}

The following theorem establishes that the candidate $(\MM,\T(\MM,\A,\mu),\partial)$ is indeed a braided associative crossed module. 
\begin{theorem}\label{P:center}
	Let $(\MM,\A,\mu)$	be an  associative crossed module.	
		 The linear map
		\[
		 \partial \colon \MM \longrightarrow \T(\MM,\A,\mu),
		\qquad \partial(m) = (\LL_m,\RR_m,\mu (m)),
		\]
	together with the action of $\T(\MM,\A,\mu)$ on $\MM$ given by
		\begin{equation} \label{action_center}
			 {}^{(\lam,\rho,a)} m \coloneqq   {}^a m, \qquad m^{(\lam,\rho,a)}  \coloneqq m^a,
		\end{equation}	 
		 and together with the bilinear map
		\[
		\{-,-\}\colon
		\T(\MM,\A,\mu) \times  \T(\MM,\A,\mu) \longrightarrow \MM
		\]
			\begin{equation}\label{braiding_center}
		\{(\lam,\rho,a),(\lam',\rho',b)\}
		\coloneqq
		\lam(b) -\lam' (a),
			\end{equation}
	forms a braided  associative crossed module.
\end{theorem}
\begin{proof}
	We begin by proving that $\partial$ is an associative morphism:
	\[
		\partial(m)\partial(n)=(\LL_m,\RR_m,\mu (m))(\LL_n,\RR_n,\mu (n))=\partial\LL_m(\mu(n))=\partial(m^{\mu(n)})\stackrel{\eqref{CM2}}{=}\partial(mn).
	\]
	
	Moreover, $\partial$ is $\T(\MM,\A,\mu)$-equivariant. Indeed:
	\[
		(\lam,\rho,a)\partial(m)=(\lam,\rho,a)(\LL_m,\RR_m,\mu(m))=\partial(\lambda\mu(m))\stackrel{\eqref{Z3}}{=}\partial({}^{a}m)=\partial({}^{(\lambda,\rho,a)}m).
	\]
	By similar arguments, one can show that $\partial(m^{(\lambda,\rho,a)})=\partial(m)(\lambda,\rho,a)$.
	
	The Peiffer identity is obtained as follows:
	\[
		{}^{\partial(m)}n={}^{(\LL_m,\RR_m,\mu(m))}n={}^{\mu(m)}n\stackrel{\eqref{CM2}\text{ for $(\MM,\A,\mu)$}}{=}mn.
	\]
	Analogously, it is also obtained that $m^{\partial(n)}=mn$. Thus, $(\MM,\T(\MM,\A,\mu),\partial)$ is an associative crossed module.
	
	Finally, let us prove that $\{-,-\}$ is a braiding. For $x=(\lam,\rho,a)$ and $y=(\lam',\rho',b)$ in $\T(\MM,\A,\mu)$, consider the bilinear map
	\[
	\langle x,y\rangle\coloneqq\lam(b).
	\]
	It is a multiplicative lifting on the crossed module $(\MM,\T(\MM,\A,\mu),\partial)$. Indeed, let $z=(\lam'',\rho'',c)\in\T(\MM,\A,\mu)$ and $m\in\MM$. Then:
	\begin{itemize}
		\item[\eqref{ML1}:] $\partial\langle x,y\rangle=\partial\lam(b)=xy$, by \eqref{eq:prod}.
		\item[\eqref{ML2}:] $\langle x,\partial(m)\rangle=\lam\mu(m)\stackrel{\eqref{Z3}}{=}{}^{a}m={}^{x}m$.
		\item[\eqref{ML3}:] $\langle\partial(m),x\rangle=\LL_m(a)=m^{a}=m^{x}$.
		\item[\eqref{ML4}:] Since $xy=\partial\lam(b)=(\LL_{\lam(b)},\RR_{\lam(b)},ab)$ and $yz=\partial\lam'(c)=(\LL_{\lam'(c)},\RR_{\lam'(c)},bc)$, we get $\langle xy,z\rangle=\LL_{\lam(b)}(c)=\lam(b)^{c}\stackrel{\eqref{Mlemu}}{=}\lam(bc)=\langle x,yz\rangle$.
	\end{itemize}
	Since $\{x,y\}=\langle x,y\rangle-\langle y,x\rangle$, Proposition~\ref{P:ml-braiding} shows that $\{-,-\}$ is a braiding on $(\MM,\T(\MM,\A,\mu),\partial)$.
\end{proof}

The lifting used in the previous proof has a companion, obtained by using the right multiplications instead of the left ones.

\begin{proposition}\label{P:liftings}
	Let $(\MM,\A,\mu)$ be an associative crossed module. For $x=(\lam,\rho,a)$ and $y=(\lam',\rho',b)$ in $\T(\MM,\A,\mu)$, the bilinear maps
	\[
	\langle x,y\rangle\coloneqq\lam(b)\qquad\text{and}\qquad\langle\!\langle x,y\rangle\!\rangle\coloneqq\rho'(a)
	\]
	are multiplicative liftings on the crossed module $(\MM,\T(\MM,\A,\mu),\partial)$ of Theorem~\ref{P:center}. Moreover, $\langle x,y\rangle-\langle\!\langle x,y\rangle\!\rangle\in\Ker(\partial)$ for all $x,y\in\T(\MM,\A,\mu)$.
\end{proposition}
\begin{proof}
	The first assertion for $\langle-,-\rangle$ is proved in Theorem~\ref{P:center}. For $\langle\!\langle-,-\rangle\!\rangle$, let $z=(\lam'',\rho'',c)\in\T(\MM,\A,\mu)$ and $m\in\MM$. Axiom \eqref{ML1} holds because $\partial\rho'(a)=xy$ by \eqref{eq:prod}. Axioms \eqref{ML2} and \eqref{ML3} read $\langle\!\langle x,\partial(m)\rangle\!\rangle=\RR_m(a)={}^{a}m$ and $\langle\!\langle\partial(m),x\rangle\!\rangle=\rho\mu(m)\stackrel{\eqref{Z4}}{=}m^{a}$.
	 For \eqref{ML4}, we have $\langle\!\langle xy,z\rangle\!\rangle=\rho''(ab)$, while $yz=\partial\rho''(b)=(\LL_{\rho''(b)},\RR_{\rho''(b)},bc)$ by Lemma~\ref{L:comp},
	 so that $\langle\!\langle x,yz\rangle\!\rangle=\RR_{\rho''(b)}(a)={}^{a}\rho''(b)\stackrel{\eqref{Rrimu}}{=}\rho''(ab)$. Finally, $\partial\langle x,y\rangle=\partial\lam(b)=\partial\rho'(a)=\partial\langle\!\langle x,y\rangle\!\rangle$ by Lemma~\ref{L:comp}.
\end{proof}
\begin{remark}\label{R:diff}
	Under the assumptions of Theorem~\ref{P:center}, the lifting  $\langle\!\langle-,-\rangle\!\rangle$ of Proposition~\ref{P:liftings} induces, by Proposition~\ref{P:ml-braiding}, a second braiding on $(\MM,\T(\MM,\A,\mu),\partial)$, namely
	\[\lBrace (\lambda,\rho,a),(\lambda',\rho',b) \rBrace \coloneqq \rho'(a)-\rho(b).\]

	By Lemma~\ref{L:comp}, $\partial\lam(b)=\partial\rho'(a)$ and $\partial\lam'(a)=\partial\rho(b)$, so the difference
	\[
	\{(\lambda,\rho,a),(\lambda',\rho',b)\}-\lBrace(\lambda,\rho,a),(\lambda',\rho',b)\rBrace=(\lam+\rho)(b)-(\lam'+\rho')(a)
	\]
	lies in $\Ker(\partial)=\{m\in\Ker(\mu)\mid {}^am=0=m^a \ \text{for all} \ a\in\A\}$. In particular, both braidings coincide whenever $\mu$ is injective.

	However, they are different in general. Let $\A$ be the two-dimensional abelian associative algebra with basis $\{e,f\}$, and let $\MM$ be the two-dimensional associative algebra with basis $\{u,v\}$ whose only non-zero product of basis elements is $u^2=v$. The action of $\A$ on $\MM$ given by ${}^eu=v=u^e$, all the other actions of basis elements being zero, together with the linear map $\mu\colon\MM\to\A$, $\mu(u)=e$, $\mu(v)=0$, defines an associative crossed module $(\MM,\A,\mu)$. 
	
	Let $\tilde{\lam} \colon\A\to\MM$ be the linear map given by $\tilde{\lam}(e)=0$ and $\tilde{\lam}(f)=v$. Then $(\tilde{\lam},0,0)\in\T(\MM,\A,\mu)$, and moreover, $(0,0,f)\in\T(\MM,\A,\mu)$, since $f$ annihilates $\A$ and acts trivially on $\MM$. Then
	\[
	\{(\tilde{\lam},0,0),(0,0,f)\}=\tilde{\lam}(f)=v\neq0=\lBrace(\tilde{\lam},0,0),(0,0,f)\rBrace.
	\]
\end{remark}

\begin{definition}\label{D:center}
The braided  associative crossed module 
\[
\ZC(\MM,\A,\mu)
\coloneqq
\bigl(\MM,\T(\MM,\A,\mu),\partial\bigr),
\]
with the braiding \eqref{braiding_center}, is called the \emph{center} of the crossed module $(\MM,\A,\mu)$.
\end{definition}

\begin{remark}\label{R:center-lifting}
	By Proposition~\ref{P:single}, the center is determined by the linear map $\partial$ together with the multiplicative lifting $\langle-,-\rangle$ of Proposition~\ref{P:liftings}. Indeed, for $x,y\in\T(\MM,\A,\mu)$ and $m,n\in\MM$, the product of $\T(\MM,\A,\mu)$, the product of $\MM$ and the action of $\T(\MM,\A,\mu)$ on $\MM$ are recovered as
	\[
	xy=\partial\langle x,y\rangle,\qquad mn=\langle\partial(m),\partial(n)\rangle,\qquad {}^{x}m=\langle x,\partial(m)\rangle,\qquad m^{x}=\langle\partial(m),x\rangle,
	\]
	and the braiding \eqref{braiding_center} is $\{x,y\}=\langle x,y\rangle-\langle y,x\rangle$. For $x=(\lam,\rho,a)$, $y=(\lam',\rho',b)$ and $z=(\lam'',\rho'',c)$, the single axiom \eqref{E:pas} reads $\lam(b)^c=\lam\mu\lam'(c)$, which holds by \eqref{Mlemu} and \eqref{Z1}. Hence Propositions~\ref{P:single} and~\ref{P:ml-braiding} give an alternative proof of Proposition~\ref{prop:algebra} and Theorem~\ref{P:center}.

	In particular, the center of any associative crossed module admits a multiplicative lifting, although the crossed module itself need not (see Remark~\ref{R:noML}).
\end{remark}

\begin{examples}\label{E:centers} \hfill
	\begin{enumerate}
		\item Let $\II$ be an ideal of $\A$ and consider the crossed module $(\II,\A,\inc)$ of Example~\ref{E:CM}~(1). Since $\inc$ is injective, \eqref{Z1} and \eqref{Z2} force $\lam(b)=ab$ and $\rho(b)=ba$ for every $(\lam,\rho,a)\in\T(\II,\A,\inc)$ and $b\in\A$; in particular, $a$ belongs to the ideal
		\[
		(\II:\A)\coloneqq\{a\in\A \mid a\A+\A a\subseteq\II\}
		\]
		of $\A$, which contains $\II$. Conversely, for any $a\in(\II:\A)$, the associativity of $\A$ shows that the triple $(\lam_a,\rho_a,a)$, with $\lam_a(b)=ab$ and $\rho_a(b)=ba$, belongs to $\T(\II,\A,\inc)$. Hence $\omega$ is an isomorphism of associative algebras from $\T(\II,\A,\inc)$ onto $(\II:\A)$ that maps $\partial(m)$ to $m$, and
		\[
		\ZC(\II,\A,\inc)\cong\bigl(\II,(\II:\A),\inc\bigr),
		\]
		with the braiding \eqref{braiding_center} given by $\{a,b\}=ab-ba$. In particular, $\ZC(\A,\A,\id)\cong(\A,\A,\id)$ with the braiding of Example~\ref{E:braidings}~(1), and $\ZC(0,\A,0)\cong(0,\ZZ(\A),0)$, since $(0:\A)=\ZZ(\A)$. For the algebra $\A$ of Example~\ref{E:CM}~(3), we have $(\K f:\A)=\K f$, so $\ZC(\K f,\A,\inc)\cong(\K f,\K f,\id)$, while $\ZC(0,\A,0)=0$ because $\ZZ(\A)=0$.

		\item Let $\MM$ be an $\A$-bimodule and consider the crossed module $(\MM,\A,0)$ of Example~\ref{E:CM}~(2). For $(\lam,\rho,a)\in\T(\MM,\A,0)$, conditions \eqref{Z1}--\eqref{Z4} say exactly that $ab=0=ba$ for all $b\in\A$ and ${}^am=0=m^a$ for all $m\in\MM$, that is, $a\in\ZZ_{\MM}(\A)$ (see Section~\ref{S:Huq}), while $(\lam,\rho)\in\Mul(\A,\MM)$ is arbitrary. Hence
		\[
		\T(\MM,\A,0)=\{(\lam,\rho,a)\mid (\lam,\rho)\in\Mul(\A,\MM),\ a\in\ZZ_{\MM}(\A)\},
		\]
		so that $\T(\MM,\A,0)\cong\Mul(\A,\MM)\oplus\ZZ_{\MM}(\A)$ as vector spaces, and the map $\partial$ is given by $\partial(m)=(\LL_m,\RR_m,0)$. Moreover, the product in $\T(\MM,\A,0)$ is trivial: for $(\lam,\rho,a),(\lam',\rho',b)\in\T(\MM,\A,0)$ and $c\in\A$, we have $\lam(b)^c=\lam(bc)=0$ and ${}^c\lam(b)=\rho(c)^b=0$, so $\partial\lam(b)=0$. The action of $\T(\MM,\A,0)$ on $\MM$ is also trivial, and the braiding \eqref{braiding_center} is given by
		\[
		\{(\lam,\rho,a),(\lam',\rho',b)\}=\lam(b)-\lam'(a).
		\]

		\item Consider the crossed module $(\MM,\A,\mu)$ of Example~\ref{E:CM}~(3), and let $(\lam,\rho,a)\in\T(\MM,\A,\mu)$ with $a=k_1 e+k_2 f$. By \eqref{Z1}, $k_1 e=ae=\mu\lam(e)\in\Im(\mu)=\K f$, so $k_1=0$. Then \eqref{Z1} and \eqref{Z2} yield $\lam(\A)\subseteq\Ker(\mu)=\K v$, $\rho(e)\in k_2 u+\K v$ and $\rho(f)\in\K v$, and \eqref{Z3}, \eqref{Z4} with $m=u$ give $\lam(f)={}^{k_2 f}u=0$ and $\rho(f)=u^{k_2 f}=0$. Finally, since ${}^ev=0=v^e$, \eqref{Mlemu} and \eqref{Rrimu} yield $\lam(e)=\lam(e)^e=0$ and $\rho(e)={}^e\rho(e)=k_2 u$. Conversely, it is straightforward to check that these triples belong to $\T(\MM,\A,\mu)$. Therefore,
		\[
		\T(\MM,\A,\mu)=\{(0,\rho_{k_2},k_2 f)\mid k_2\in\K\}, \qquad \text{where} \ \ \rho_{k_2}(e)=k_2 u, \ \ \rho_{k_2}(f)=0.
		\]
		Since $\partial(u)=(\LL_u,\RR_u,f)=(0,\rho_1,f)$ and $\partial(v)=0$, the map $\partial$ is surjective with kernel $\K v$. The product in $\T(\MM,\A,\mu)$ is trivial, and $\omega$ is an isomorphism from $\T(\MM,\A,\mu)$ onto the abelian subalgebra $\K f$ of $\A$, which acts trivially on $\MM$. Consequently,
		\[
		\ZC(\MM,\A,\mu)\cong(\MM,\K f,\mu),
		\]
		and the braiding \eqref{braiding_center} is trivial, since $\lam=0$ for all $(\lam,\rho,a)\in\T(\MM,\A,\mu)$. In particular, $\pi_1\big(\ZC(\MM,\A,\mu)\big)=\K v=\pi_1(\MM,\A,\mu)$, whereas $\pi_0\big(\ZC(\MM,\A,\mu)\big)=0$ and $\pi_0(\MM,\A,\mu)=\A/\K f\cong\K$.	
	\end{enumerate}
\end{examples}

Besides, the associative morphism $\omega$ introduced in Lemma \ref{L:omega_mor} becomes an associative crossed module with a suitable action.
\begin{proposition} \label{P:omega}
	Let $(\MM,\A,\mu)$	be an   associative crossed module.  The map 
\[	\omega \colon \T(\MM,\A,\mu) \rightarrow   \A, \quad   \omega(\lam,\rho,a) \coloneqq a,\]
 constitutes an  associative crossed module, where the action of $\A$ is defined by 
		\begin{equation}\label{action_L}
			{}^b  (\lam,\rho,a) \coloneqq  (\LL_{\rho(b)},\RR_{\rho(b)}, ba ) =\partial \rho (b), \ (\lam,\rho,a)^b  \coloneqq   (\LL_{\lam(b)},\RR_{\lam(b)},  ab)=\partial \lam (b). 
		\end{equation}	
\end{proposition}
\begin{proof}
	Let us prove that the maps \eqref{action_L} define an action of $\A$ on $\T(\MM,\A,\mu)$:
	\begin{itemize}
		\item[\eqref{A1}:] ${}^{b}({}^{c}(\lambda,\rho,a))={}^{b}\partial\rho(c)=\partial\RR_{\rho(c)}(b)=\partial({}^{b}\rho(c))
		\stackrel{\eqref{Rrimu}}{=}\partial\rho(bc)={}^{bc}(\lambda,\rho,a)$.
		\item[\eqref{A2}:] $((\lambda,\rho,a)^{b})^{c}=\partial\lambda(b)^{c}=\partial\LL_{\lambda(b)}(c)=\partial(\lambda(b)^{c})
		\stackrel{\eqref{Mlemu}}{=}\partial\lambda(bc)=(\lambda,\rho,a)^{bc}$.
		\item[\eqref{A3}:]
		\begin{align*}
			({}^{c}(\lambda,\rho,a))(\lambda',\rho',b)&=\partial\rho(c)(\lambda',\rho',b)=\partial\LL_{\rho(c)}(b)
			=\partial\rho(c)^{b}\stackrel{\eqref{Mbimul}}{=}\partial\,{}^{c}\lambda(b)\\
			&=\partial\RR_{\lambda(b)}(c)={}^{c}\partial\lambda(b)={}^{c}\big((\lambda,\rho,a)(\lambda',\rho',b)\big).
		\end{align*}
		\item[\eqref{A4}:] 
		\begin{align*}
			(\lambda,\rho,a)((\lambda',\rho',b)^{c})&=(\lambda,\rho,a)\partial\lambda'(c)=\partial\lambda(bc)
			\stackrel{\eqref{Mlemu}}{=}\partial\lambda(b)^{c}=\partial\LL_{\lambda(b)}(c)\\
			&=(\partial\lambda(b))^{c}=\big((\lambda,\rho,a)(\lambda',\rho',b)\big)^{c}.
		\end{align*}
		\item[\eqref{A5}:] 
		\begin{align*}
			({}^{b}(\lambda,\rho,a))^{c}&=(\partial\rho(b))^{c}=\partial\LL_{\rho(b)}(c)=\partial\rho(b)^{c}
			\stackrel{\eqref{Mbimul}}{=}\partial\,{}^{b}\lambda(c)=\partial\RR_{\lambda(c)}(b)\\
		   &={}^{b}\partial\lambda(c)={}^{b}((\lambda,\rho,a)^{c}).
		\end{align*}
		\item[\eqref{A6}:]
		\begin{align*}
			((\lambda,\rho,a)^{c})(\lambda',\rho',b)&=\partial\lambda(c)(\lambda',\rho',b)=\partial\LL_{\lambda(c)}(b)
			=\partial\lambda(c)^{b}\stackrel{\eqref{Mlemu}}{=}\partial\lambda(cb)\\
			&=(\lambda,\rho,a)\partial\rho'(c)=(\lambda,\rho,a)({}^{c}(\lambda',\rho',b)).
		\end{align*}
	\end{itemize}
	
	Therefore, $(\T(\MM,\A,\mu),\A,\omega)$ is an associative crossed module. Indeed, $\omega$ is $\A$-equivariant:
	\begin{align*}
		&\omega({}^{b}(\lambda,\rho,a))=\omega(\partial\rho(b))=ba=b\omega(\lambda,\rho,a),\\
		&\omega((\lambda,\rho,a)^{b})=\omega(\partial\lambda(b))=\omega(ab)=\omega(\lambda,\rho,a)b,
	\end{align*}
	and the Peiffer identity also holds:
	\begin{align*}
		&{}^{\omega(\lambda,\rho,a)}(\lambda',\rho',b)={}^{a}(\lambda',\rho',b)=\partial\rho'(a)\stackrel{\text{Lemma \ref{L:comp}}}{=}(\lambda,\rho,a)(\lambda',\rho',b),\\
		&(\lambda,\rho,a)^{\omega(\lambda',\rho',b)}=(\lambda,\rho,a)^{b}=\partial\lambda(b)=(\lambda,\rho,a)(\lambda',\rho',b).\qedhere
	\end{align*}
\end{proof}

As a direct consequence of the definition \eqref{action_center} of the action and of the equality $\omega\partial=\mu$, we have the following corollary.
\begin{corollary}\label{C:forget}
The pair $(\id,\omega) \colon \ZC(\MM,\A,\mu) \to  (\MM,\A,\mu) $ is a  morphism of associative crossed modules.
\end{corollary}
\begin{corollary}[Homotopy fiber sequence] \label{C:fiber}
	There is an exact sequence of homotopy groups
	\vspace{-0.4cm}
		\[
	\begin{tikzpicture}
		\matrix (m) [
		matrix of math nodes,
		row sep=1.5em,
		column sep=0.01cm,
		text height=2.8ex, text depth=1.5ex
		]
		{ 0 & \pi_1\big(\ZC(\MM,\A,\mu) \big)  & \pi_1(\MM,\A,\mu)  & \pi_1\big(\T(\MM,\A,\mu),\A,\omega \big) \\
			\pi_0\big(\ZC(\MM,\A,\mu) \big) &\pi_0(\MM,\A,\mu) & \pi_0\big(\T(\MM,\A,\mu),\A,\omega \big) & 0. \\
		};
		\path[overlay,->, font=\scriptsize,>=angle 90]
		(m-1-1) edge (m-1-2)
		(m-1-2) edge (m-1-3)
		(m-1-3) edge (m-1-4)
		(m-1-4) edge[out=355,in=175]  (m-2-1)
		(m-2-1) edge (m-2-2)
		(m-2-2) edge (m-2-3)
		(m-2-3) edge (m-2-4);
	\end{tikzpicture}
		\]
\end{corollary}
\begin{proof}
The proof follows from Corollary~\ref{C:forget},  the equality $\omega \partial = \mu$, and the application of the kernel–cokernel lemma to the composition  $\omega \partial$.
\end{proof}
\begin{proposition}\label{P:square} The square
\[	
	\begin{tikzcd}
		\MM \arrow[d,equal,"\id",swap] \arrow[r, "\partial"]  & \T(\MM,\A,\mu) \arrow[d,"\omega"] \\
		\MM \arrow[r,"\mu",swap]                        & \A           
	\end{tikzcd} 
\]
together with
\begin{align*}
	h_1 \colon \MM \times \T(\MM,\A,\mu) \xrightarrow{ \quad } \MM, \ h_1(m,(\lam,\rho,a)) \coloneqq m^{(\lam,\rho,a)} = m^a,\\
	h_2 \colon  \T(\MM,\A,\mu) \times \MM  \xrightarrow{ \quad } \MM, \ h_2((\lam,\rho,a),m) \coloneqq	 {}^{(\lam,\rho,a)} m= {}^a m.
\end{align*}
	is a crossed square of associative algebras.
\end{proposition}
\begin{proof}
	We will show that the previous square satisfies the conditions in Definition~\ref{square}.
    The square is commutative because $\mu= \omega \partial$, $\A$ is acting on $\MM$ and $\T(\MM,\A,\mu)$ (see \eqref{action_L}), $\MM$ acts on $\MM$ by multiplication, $\T(\MM,\A,\mu)$ acts on $\MM$ (recall \eqref{action_center}) and the maps $h_1,h_2$ are bilinear. Furthermore, we obtain the following identities, with all omitted equalities being immediate.
    
    $\bullet$ \eqref{SC1}: It is a consequence of Theorem~\ref{P:center}, Proposition~\ref{P:omega} and Corollary~\ref{C:forget}.

$\bullet$ \eqref{SC2}: Consider the action of $\MM$ on $\T(\MM,\A,\mu)$ induced by $\mu$, that is, 
\[{}^{m}(\lambda,\rho,a)\coloneqq{}^{\mu(m)}(\lambda,\rho,a) \quad
\text { and } \quad  (\lambda,\rho,a)^{m}\coloneqq(\lambda,\rho,a)^{\mu(m)}.\] 
    We have that 
		\[\partial(h_1(m,(\lambda,\rho,a)))=\partial(m^a)\stackrel{\eqref{Z4}}{=}\partial(\rho\mu(m))\stackrel{\eqref{action_L}}{=}{}^{\mu(m)}(\lambda,\rho,a)={}^{m}(\lambda,\rho,a).\]
   	Similarly, $\partial(h_2((\lambda,\rho,a),m))=(\lambda,\rho,a)^{m}$.
   	
$\bullet$ \eqref{SC3}: $h_1(m,\partial(m'))=m^{(\LL_{m'},\RR_{m'}, \mu(m'))}\stackrel{\eqref{action_center}}{=}m^{\mu(m')}\stackrel{\eqref{CM2}}{=}mm'$,  
 and an analogous argument yields $h_2(\partial(m),m')=mm'$.
	
$\bullet$ \eqref{SC4}: This axiom is derived from:
    \begin{align*}
		 h_1(mm', (\lambda, \rho, a))&= (mm')^{(\lambda, \rho, a)} \stackrel{\eqref{action_center}}{=} (mm')^{a}\stackrel{\eqref{A4}}{=} m(m'^{a})=mh_{1}(m',(\lambda,\rho,a)),\\
		 h_1(m, (\lambda, \rho, a)(\lambda', \rho', b)) &=m^{(\lambda, \rho, a)(\lambda', \rho', b)}=m^{\partial\lambda(b)}\stackrel{\eqref{action_center}}{=}m^{ab}\stackrel{\eqref{A2}}{=}(m^{a})^{b}\\
		 &\stackrel{\eqref{action_center}}{=} \left(m^{(\lambda, \rho, a)}\right)^{(\lambda', \rho', b)} =h_1(m, (\lambda, \rho, a))^{(\lambda', \rho', b)},\\
        h_2( (\lambda, \rho, a),m)^{(\lambda', \rho', b)} &=({}^{(\lambda, \rho, a)}m)^{(\lambda', \rho', b)}\stackrel{\eqref{action_center}}{=}({}^{a}m)^{b}\stackrel{\eqref{A5}}{=}{}^{a}(m^{b})\\
        & \stackrel{\eqref{action_center}}{=} {}^{(\lambda, \rho, a)}h_1( m,(\lambda', \rho', b)),  \\
        h_1(m, (\lambda, \rho, a))^{m'}&=  (m^{(\lambda, \rho, a)})m' 
        \stackrel{\eqref{action_center}}{=}(m^{a})m'  \stackrel{\eqref{A6}}{=} m({}^{a}m')\\
        &=mh_2( (\lambda, \rho, a),m')
        ={}^{m}h_2( (\lambda, \rho, a),m'),
	\end{align*}
	where the corresponding first two equalities above for $h_{2}$ follow readily by analogous arguments.

$\bullet$ \eqref{SC5}: The remaining axiom is obtained as follows:
\begin{align*}
		h_1(m, (\lambda, \rho, a)^{b}) &\stackrel{\eqref{action_L}}{=} m^{\partial\lambda(b)}\stackrel{\eqref{action_center}} = m^{ab}\stackrel{\eqref{A2}}{=} (m^{a})^{b}\stackrel{\eqref{action_center}} = (m^{(\lambda, \rho, a)})^{b} = h_1(m, (\lambda, \rho, a))^{b}, \\
	h_1({}^{b}m, (\lambda, \rho, a)) &= ({}^{b}m)^{(\lambda, \rho, a)} \stackrel{\eqref{action_center}}{=}({}^{b}m)^{a}\stackrel{\eqref{A5}}{=}{}^{b}(m^{a})= {}^{b}(m^{(\lambda, \rho, a)}) = {}^{b}h_1(m, (\lambda, \rho, a)),
\end{align*}
where the other identities follow analogously from the properties of the action \eqref{A1}--\eqref{A6}.
\end{proof}


\section{Associative and Lie centers: a comparison}\label{S:comp}

In this section, we compare our center with Pirashvili's center of Lie
crossed modules~\cite{Pim2} via the (braided) Liezation functors~\cite{FFL1}. We obtain a canonical
comparison morphism, which preserves braidings under an injectivity condition.

Throughout this section we assume that $\operatorname{char}(\K)\neq 2$, as in \cite{Pim2}. An associative crossed module $(\MM,\A,\mu)$ yields a \emph{braided crossed module of Lie algebras} (recall this notion in \cite[Definition 4.4]{FFL1}) in two ways:
	\begin{itemize}
		\item[(i)] Apply the \emph{Liezation} first \cite[Lemma 1.1]{FFL1}, and then take the center constructed by M. Pirashvili in \cite[Proposition 3.2]{Pim2} for a crossed module of Lie algebras. In this case we obtain
		\[\ZC_{\text{Lie}}(\MM_{\text{Lie}},\A_{\text{Lie}},\mu)=(\MM_{\text{Lie}},\T_{\text{Lie}}(\MM_{\text{Lie}},\A_{\text{Lie}},\mu),\delta)\]
		where
		$\MM_{\text{Lie}}=(\MM,[m,n])$, $\A_{\text{Lie}}=(\A,[a,b])$, and $\A_{\text{Lie}}$ acts on $\MM_{\text{Lie}}$ by $[a,m]\coloneqq{}^am-m^a$, whereas
		\[\T_{\text{Lie}}(\MM_{\text{Lie}},\A_{\text{Lie}},\mu)=\{(\xi,x)\mid x\in \A_{\text{Lie}},\,\xi\colon \A_{\text{Lie}}\to \MM_{\text{Lie}} \text{ is a derivation satisfying }\eqref{ZL1}\text{--}\eqref{ZL2}\},\]
		\begin{align}
			\mu\xi(a)&=[x,a],\tag{ZL1}\label{ZL1}\\
			\xi\mu(m)&=[x,m],\tag{ZL2}\label{ZL2}		
		\end{align}
		$\delta(m)\coloneqq(\Ad_{m},\mu(m))$ with $\Ad_{m}(a)\coloneqq- [a,m]$ for all $a\in\A,\,m\in \MM$.
		
		Note  that $\T_{\text{Lie}}(\MM_{\text{Lie}},\A_{\text{Lie}},\mu)$ is a Lie algebra with the bracket given by $[(\xi,x),(\eta,y)]=\delta\xi(y)=(\Ad_{\xi(y)},[x,y])$ (see \cite[Corollary~3.4]{Pim2}),  $\T_{\text{Lie}}(\MM_{\text{Lie}},\A_{\text{Lie}},\mu)$ acts on $\MM_{\text{Lie}}$ via $[(\xi,x),m]=[x,m]$,
		 and the braiding is given by
		$\ldblbrace(\xi,x),(\eta,y)\rdblbrace=\xi(y)$ (see \cite[Proposition~3.2]{Pim2}).

	Pirashvili considers braided crossed modules of Lie algebras as defined in \cite[Proposition~6.20]{CC}. Direct verification shows that the bracket
	 $\{(\xi,x),(\eta,y)\}=\xi(y)$ satisfies the axioms outlined in \cite[Definition~4.4]{FFL1}. Therefore, both constructions presented in this section yield braided crossed modules of Lie algebras in the same sense.

		\item[(ii)] Take the center in the sense of Definition \ref{D:center}, and then apply the \emph{braided Liezation} \cite[Proposition 4.3]{FFL1}, which requires $\operatorname{char}(\K)\neq 2$. Following this approach, we obtain the braided crossed module of Lie algebras
		\[(\MM_{\text{Lie}},\T(\MM,\A,\mu)_{\text{Lie}},\partial)\]
		with the braiding given by
		\[\{(\lambda,\rho,a),(\lambda',\rho',b)\}_{\text{Lie}}=\frac{1}{2}(\{(\lambda,\rho,a),(\lambda',\rho',b)\}-\{(\lambda',\rho',b),(\lambda,\rho,a)\})=\lambda(b)-\lambda'(a),\]
		i.e. in this case $\{(\lambda,\rho,a),(\lambda',\rho',b)\}_{\text{Lie}}=\{(\lambda,\rho,a),(\lambda',\rho',b)\}$.
	\end{itemize}
	
	Therefore, a natural question to tackle is to study the relationship between both constructions. Upon forgetting the braidings,
	 there is a canonical morphism between the underlying crossed modules in the Lie setting.
	
	\begin{proposition}\label{P:canmor}
		Let $(\MM,\A,\mu)$ be an associative crossed module. The pair 
		\[( \id,\Phi)\colon (\MM_{\textnormal{Lie}},\T(\MM,\A,\mu)_{\textnormal{Lie}},\partial)\rightarrow \ZC_{\textnormal{Lie}}(\MM_{\textnormal{Lie}},\A_{\textnormal{Lie}},\mu)\]
		 is a morphism of crossed modules of Lie algebras, where $\Phi(\lambda,\rho,a)=(\lambda-\rho,a)$. 
	\end{proposition}
	\begin{proof}
		We begin by showing that $\Phi$ is well-defined, i.e. $\lambda-\rho$ is a derivation and  $(\lambda-\rho,a)$ satisfies
		 \eqref{ZL1}--\eqref{ZL2} for all $(\lambda,\rho,a)\in\T(\MM,\A,\mu)_{\text{Lie}}$. Let $b,c\in\A$.
		
		$\bullet$	$\lambda-\rho$ derivation. 
			\begin{align*}
				[b,(\lambda-\rho)(c)]-[c,(\lambda-\rho)(b)]&=[b,\lambda(c)]-[b,\rho(c)]-[c,\lambda(b)]+[c,\rho(b)]\\
				&={}^b\lambda(c)-\lambda(c)^{b}-{}^b\rho(c)+\rho(c)^b\\
				&\quad-{}^c\lambda(b)+\lambda(b)^c+{}^c\rho(b)-\rho(b)^{c}\\
				&\stackrel{\eqref{Mbimul}}{=}\lambda(b)^c-\lambda(c)^b-{}^{b}\rho(c)+{}^{c}\rho(b)\\
				&\stackrel{\eqref{Mlemu},\eqref{Rrimu}}{=}(\lambda-\rho)(bc-cb)=(\lambda-\rho)[b,c].
			\end{align*}
			$\bullet$ \eqref{ZL1}. $\mu(\lambda-\rho)(b)=\mu\lambda(b)-\mu\rho(b)\stackrel{\eqref{Z1},\eqref{Z2}}{=}ab-ba=[a,b]$.
			
		\noindent	$\bullet$\eqref{ZL2}. Let $m\in \MM$. $(\lambda-\rho)\mu(m)=\lambda\mu(m)-\rho\mu(m)\stackrel{\eqref{Z3},\eqref{Z4}}{=}{}^{a}m-m^{a}=[a,m]$.

		Clearly, $\Phi$ is linear. Let us prove that it also preserves the Lie bracket: on the one hand,
		\begin{align*}
			\Phi[(\lambda,\rho,a),(\lambda',\rho',b)]&=\Phi((\lambda,\rho,a)(\lambda',\rho',b)-(\lambda',\rho',b)(\lambda,\rho,a))
			=\Phi\partial\lambda(b)-\Phi\partial\lambda'(a)\\
			&=\Phi\partial(\lambda-\rho)(b)\;\text{\footnotesize(by Lemma~\ref{L:comp}, $\partial\lambda'(a)=\partial\rho(b)$)}\\
			&=(\LL_{(\lambda-\rho)(b)}-\RR_{(\lambda-\rho)(b)},\mu(\lambda-\rho)(b))
		\end{align*}
		and, on the other hand
		\begin{align*}
			[\Phi(\lambda,\rho,a),\Phi(\lambda',\rho',b)]&=[(\lambda-\rho,a),(\lambda'-\rho',b)]\\
			&=\delta(\lambda-\rho)(b)=(\Ad_{(\lambda-\rho)(b)},\mu(\lambda-\rho)(b)).
		\end{align*}
		Thus, $\Phi$ is a Lie morphism if and only if $\LL_{(\lambda-\rho)(b)}-\RR_{(\lambda-\rho)(b)}=\Ad_{(\lambda-\rho)(b)}$. Indeed, given $x\in\A$,
		\begin{align*}
			\LL_{(\lambda-\rho)(b)}(x)-\RR_{(\lambda-\rho)(b)}(x)&=(\lambda-\rho)(b)^x-{}^{x}(\lambda-\rho)(b)\\
			&=-[x,(\lambda-\rho)(b)]=\Ad_{(\lambda-\rho)(b)}(x).
		\end{align*}
		
		So, since the compatibility between the actions is direct, it only remains to show that $\Phi\partial=\delta$, which follows by
		\begin{align*}
			\Phi\partial(m)=\Phi(\LL_m,\RR_m,\mu(m))&=(\LL_m-\RR_m,\mu(m))\\
			&=(-[-,m],\mu(m))=(\Ad_{m},\mu(m))=\delta(m).\qedhere
		\end{align*}
	\end{proof}

However, the canonical morphism $( \id,\Phi)$ is not necessarily a morphism of braided crossed modules of Lie algebras. For this to hold,
 it must preserve the braiding, that is, for all $(\lambda,\rho,a),(\lambda',\rho',b)\in\T(\MM,\A,\mu)$,
\begin{align*}
	&\{(\lambda,\rho,a),(\lambda',\rho',b)\}_{\text{Lie}}=\{\Phi(\lambda,\rho,a),\Phi(\lambda',\rho',b)\}\\
	\iff&\lambda(b)-\lambda'(a)=\{(\lambda-\rho,a),(\lambda'-\rho',b)\}=(\lambda-\rho)(b)\iff\lambda'(a)=\rho(b).
\end{align*}
In terms of the multiplicative liftings of Proposition~\ref{P:liftings}, this gives the following characterization.
\begin{proposition}\label{P:braided-iff}
	Let $(\MM,\A,\mu)$ be an associative crossed module. The morphism $(\id,\Phi)$ in Proposition \ref{P:canmor} is braided if and only if
	 the two multiplicative liftings $\langle-,-\rangle$ and $\langle\!\langle-,-\rangle\!\rangle$ of Proposition~\ref{P:liftings} coincide,
	  that is, $\lam(b)=\rho'(a)$ for all $(\lambda,\rho,a),(\lambda',\rho',b)\in\T(\MM,\A,\mu)$. In particular, $(\id,\Phi)$ is braided whenever the map $\partial$ is injective.
\end{proposition}
\begin{proof}
	By the computation above, $(\id,\Phi)$ is braided if and only if $\lam'(a)=\rho(b)$ for all $x=(\lambda,\rho,a)$ and $y=(\lambda',\rho',b)$ in $\T(\MM,\A,\mu)$,
	 that is, $\langle y,x\rangle=\langle\!\langle y,x\rangle\!\rangle$ for all $x,y$. If $\partial$ is injective, this holds because
	  $\langle y,x\rangle-\langle\!\langle y,x\rangle\!\rangle\in\Ker(\partial)$ by Proposition~\ref{P:liftings}.
\end{proof}
\begin{example}
	Let $(\MM,\A,\mu)$ be the associative crossed module of Remark~\ref{R:diff}, and let $\tilde{\lam}\colon\A\to\MM$ be the linear map given there,
	 that is, $\tilde{\lam}(e)=0$ and $\tilde{\lam}(f)=v$. We show that the morphism $(\id,\Phi)$ of Proposition~\ref{P:canmor} is neither an isomorphism nor braided.
	  Note that $\partial$ is not injective, since $\partial(v)=0$, so the sufficient condition in Proposition~\ref{P:braided-iff} does not hold.
	
	Let $(\lam,\rho,a)\in\T(\MM,\A,\mu)$ with $a=k_1 e+k_2 f$, $k_1,k_2\in\K$. Since the product of $\A$ is trivial, \eqref{Z1} and \eqref{Z2} give $\lam(\A),\rho(\A)\subseteq\Ker(\mu)=\K v$,
	 while \eqref{Z3} and \eqref{Z4} give $\lam(e)=\lam\mu(u)={}^au=k_1 v$ and $\rho(e)=\rho\mu(u)=u^a=k_1 v$. Conversely, any such triple satisfies \eqref{Z1}--\eqref{Z4}, 
	 \eqref{Mlemu}, \eqref{Rrimu} and \eqref{Mbimul}, because $\A$ has trivial product and $v$ is annihilated by $\A$. Hence $\T(\MM,\A,\mu)$ consists of the triples 
	 $(\lam,\rho,k_1 e+k_2 f)$, with $k_1,k_2\in\K$, such that
	\[
	\lam(e)=\rho(e)=k_1 v \quad\text{and}\quad \lam(f),\rho(f)\in\K v,
	\]
	so it has dimension $4$.
	
	On the other hand, $\MM_{\text{Lie}}$ and $\A_{\text{Lie}}$ are abelian Lie algebras, and the action of $\A_{\text{Lie}}$ on $\MM_{\text{Lie}}$ is trivial,
	 since ${}^am=m^a$ for all $a\in\A$ and $m\in\MM$. Equations \eqref{ZL1} and \eqref{ZL2} reduce to $\mu\xi=0$ and $\xi\mu=0$. Since $\mu(\MM)=\K e$, we get
	\[
	\T_{\text{Lie}}(\MM_{\text{Lie}},\A_{\text{Lie}},\mu)=\{(c,\xi)\mid c\in\A,\ \xi\colon\A\to\MM \text{ linear},\ \xi(e)=0,\ \xi(f)\in\K v\},
	\]
	which has dimension $3$.
	
	Now $\Phi(\lam,\rho,k_1 e+k_2 f)=(\lam-\rho,k_1 e+k_2 f)$, where $(\lam-\rho)(e)=0$ and $(\lam-\rho)(f)=\lam(f)-\rho(f)$.
	 Therefore $\Phi$ is surjective, but it is not injective: its kernel is spanned by the non-zero element $(\tilde{\lam},\tilde{\lam},0)$.
	  In particular, $(\id,\Phi)$ is not an isomorphism.
	
	Finally, $(\tilde{\lam},0,0)$ and $(0,0,f)$ belong to $\T(\MM,\A,\mu)$, and
	\begin{align*}
		\{(0,0,f),(\tilde{\lam},0,0)\}_{\text{Lie}}&=0-\tilde{\lam}(f)=-v,\\
		\ldblbrace \Phi(0,0,f),\Phi(\tilde{\lam},0,0) \rdblbrace &=\ldblbrace (0,f),(\tilde{\lam},0) \rdblbrace=0.
	\end{align*}
	Since $v\neq0$, $(\id,\Phi)$ is not a morphism of braided crossed modules of Lie algebras. Equivalently, in accordance with Proposition~\ref{P:braided-iff},
	 the two multiplicative liftings differ: $\langle(\tilde{\lam},0,0),(0,0,f)\rangle=\tilde{\lam}(f)=v$, whereas $\langle\!\langle(\tilde{\lam},0,0),(0,0,f)\rangle\!\rangle=0$.
\end{example}


\section{Comparison with the Huq center}\label{S:Huq}
In this section, we construct the Huq center of an associative crossed module and compare it with the previously introduced notion of center in Definition \ref{D:center}. 

We briefly recall the categorical notion of center introduced by Huq \cite{Huq}. Let $\mathcal{C}$ be a pointed category with finite products and 
let $f\colon X\to Z$ and $g\colon Y\to Z$ be morphisms in $\mathcal{C}$. They are said to \emph{Huq-commute} if there exists a morphism 
$\varphi\colon X\times Y\to Z$ in $\mathcal{C}$ such that $\varphi\circ(\id_X,0)=f$ and $\varphi\circ(0,\id_Y)=g$, where $(s,t)$
 denotes the unique morphism into a product whose components are $s$ and $t$. The \emph{Huq center} of an object $X$, whenever it exists,
  is the largest subobject $ \ZZ^{\mathrm{Huq}}(X)\rightarrow X$ which Huq-commutes with $\id_X$.

The category of associative algebras is semi-abelian \cite{JMT}. Moreover, crossed modules in a semi-abelian category are equivalent to internal groupoids \cite{J},
 and the category of internal groupoids in a semi-abelian category is again semi-abelian \cite{BG}. Consequently, the category of crossed modules of associative algebras is semi-abelian and,
  in particular, pointed with finite products. Hence, Huq commutativity, and therefore the corresponding notion of Huq center, applies in our setting.

\begin{definition}
		Let $\A$ be an associative algebra  acting on $\MM$.  The \emph{annihilator of $\A$ on $\MM$}, denoted by
		$\Ann_{\MM}(\A)$, is defined as
		\[ \Ann_{\MM}(\A) \coloneqq \{m \in \MM \mid \LL_m=0= \RR_m \}= \{m \in \MM \mid  m^a =0= {}^a m \ \text{for all} \  a \in \A \},\]
		the \emph{annihilator of $\MM$ on $\A$}, denoted by
		$\Ann_{\A}(\MM)$, is defined as
		\[ \Ann_{\A}(\MM) \coloneqq  \{a \in \A \mid  {}^a m=0= m^a \ \text{for all} \  m \in \MM \},\]
		and the \emph{center of $\A$ on $\MM$},  denoted by
		$\ZZ_{\MM}(\A)$, is defined as
		\[ \ZZ_{\MM}(\A) \coloneqq \{a \in \ZZ(\A) \mid {}^a m=0= m^a \ \text{for all} \  m \in \MM \}= \ZZ(\A) \cap \Ann_{\A}(\MM).\]
	\end{definition} 

The following theorem provides an explicit characterization of the Huq center of an associative crossed module (see the annihilator of a crossed module in \cite{AH}; cf. also~\cite{UCA}).

\begin{theorem}
	Let $(\MM,\A,\mu)$ be an associative crossed module. The Huq center of $(\MM,\A,\mu)$ is the abelian associative crossed module
	\begin{equation} \label{HuqC}
		\ZZ^{\text{\emph{Huq}}}(\MM,\A,\mu) \coloneqq (\Ann_{\MM}(\A),\ZZ_{\MM}(\A),\mu|_{\Ann_{\MM}(\A)}).
	\end{equation}
\end{theorem}
\begin{proof}
	We first characterize the associative crossed submodules that Huq-commute with $(\MM,\A,\mu)$, and then show that \eqref{HuqC} is one of them.
	
\textbf{Step 1:} An associative crossed submodule $(\NN,\B,\nu)$ of $(\MM,\A,\mu)$ Huq-commutes with $(\MM,\A,\mu)$ if and only if $\NN\subseteq\Ann_{\MM}(\A)$ and $\B\subseteq\ZZ_{\MM}(\A)$.
	
	Recall that the product of two associative crossed modules is computed component-wise, with component-wise products and actions. If a morphism
	\[(\NN,\B,\nu)\times (\MM,\A,\mu)= (\NN \times \MM,\B \times \A,\nu \times \mu)\xrightarrow{ \ (\varphi_1,\varphi_2) \ } (\MM,\A,\mu)\]
	satisfies $\varphi_1(n,0)=n$, $\varphi_1(0,m)=m$, $\varphi_2(b,0)=b$ and $\varphi_2(0,a)=a$, then, since $\varphi_1$ and $\varphi_2$ are linear, $\varphi_1(n,m)=n+m$ and $\varphi_2(b,a)=b+a$. Hence $(\NN,\B,\nu)$ Huq-commutes with $(\MM,\A,\mu)$ if and only if the addition maps
	\[
	f\colon \NN\times \MM\rightarrow \MM,\ \ f(n,m)=n+m, \qquad g\colon \B\times \A\rightarrow \A,\ \ g(b,a)=b+a,
	\]
	form a morphism of associative crossed modules. Since $\nu=\mu|_{\NN}$, the equality $\mu f=g(\nu\times\mu)$ always holds. For $n,n'\in\NN$, $m,m'\in\MM$, $b,b'\in\B$ and $a,a'\in\A$, we have
	\begin{align*}
		g(b,a)\,g(b',a')-g\big((b,a)(b',a')\big)&=ba'+ab',\\
		f(n,m)\,f(n',m')-f\big((n,m)(n',m')\big)&=nm'+mn',\\
		{}^{g(b,a)}f(n,m)-f\big({}^{(b,a)}(n,m)\big)&={}^{b}m+{}^{a}n,\\
		f(n,m)^{g(b,a)}-f\big((n,m)^{(b,a)}\big)&=m^{b}+n^{a}.
	\end{align*}
	Setting all but two of the variables equal to zero, we see that these four expressions vanish identically if and only if
	\[
	b\A=0=\A b,\qquad n\MM=0=\MM n,\qquad {}^{b}m=0=m^{b},\qquad {}^{a}n=0=n^{a}
	\]
	for all $n\in\NN$, $m\in\MM$, $b\in\B$ and $a\in\A$. The first and third conditions say exactly that $\B\subseteq\ZZ(\A)\cap\Ann_{\A}(\MM)=\ZZ_{\MM}(\A)$, and the fourth one says that $\NN\subseteq\Ann_{\MM}(\A)$. Finally, the second condition follows from the fourth: if $n\in\Ann_{\MM}(\A)$, then, for all $m\in\MM$,
	\begin{equation}\label{Ann_prod}
		nm\stackrel{\eqref{CM2}}{=}n^{\mu(m)}=0 \qquad\text{and}\qquad mn\stackrel{\eqref{CM2}}{=}{}^{\mu(m)}n=0.
	\end{equation}
	
	\textbf{Step 2:} $(\Ann_{\MM}(\A),\ZZ_{\MM}(\A),\mu|_{\Ann_{\MM}(\A)})$ is an abelian associative crossed submodule of $(\MM,\A,\mu)$.
	
	By \eqref{Ann_prod}, $\Ann_{\MM}(\A)$ has trivial product, and so does $\ZZ_{\MM}(\A)\subseteq\ZZ(\A)$; hence both are abelian subalgebras. Moreover, $\ZZ_{\MM}(\A)$ acts trivially on $\MM$, and in particular on $\Ann_{\MM}(\A)$. It remains to show that $\mu(\Ann_{\MM}(\A))\subseteq\ZZ_{\MM}(\A)$. Let $m\in\Ann_{\MM}(\A)$, $n\in\MM$ and $a\in\A$. Then
	\begin{align*}
		a\mu(m)&\stackrel{\eqref{CM1}}{=}\mu({}^{a}m)=0, & \mu(m)a&\stackrel{\eqref{CM1}}{=}\mu(m^{a})=0,\\
		{}^{\mu(m)}n&\stackrel{\eqref{CM2}}{=}mn\stackrel{\eqref{Ann_prod}}{=}0, & n^{\mu(m)}&\stackrel{\eqref{CM2}}{=}nm\stackrel{\eqref{Ann_prod}}{=}0,
	\end{align*}
	so $\mu(m)\in\ZZ_{\MM}(\A)$.
	
	By Step 1, every associative crossed submodule that Huq-commutes with $(\MM,\A,\mu)$ is contained in \eqref{HuqC}, and, by Steps 1 and 2,
	 \eqref{HuqC} itself Huq-commutes with $(\MM,\A,\mu)$. Hence \eqref{HuqC} is the Huq center of $(\MM,\A,\mu)$.
\end{proof}

We now compare the Huq center \eqref{HuqC} of an associative crossed module with the center constructed in Definition~\ref{D:center}.

Define $\zeta  \colon \ZZ_{\MM}(\A) \to  \T(\MM,\A,\mu)$ as the injection given by $\zeta(a)=(0,0,a)$.  This map $\zeta$ is  well defined because,
 for any $a \in \ZZ_{\MM}(\A)$, it follows that $(0,0,a)\in\T(\MM,\A,\mu)$.

\begin{proposition} The map 
	\[(\inc,\zeta) \colon   \ZZ^{\textnormal{Huq}}(\MM,\A,\mu)  \xrightarrow{ \  \ } \ZC(\MM,\A,\mu)\] is an  injective morphism of associative crossed modules
	 that induces an isomorphism on $\pi_1$  and an injection on  $\pi_0$.
\end{proposition}
\begin{proof}
First, $(\inc,\zeta)$ is a morphism of associative crossed modules. Indeed, for $m\in\Ann_{\MM}(\A)$ we have $\LL_m=0=\RR_m$, so $\partial(m)=(0,0,\mu(m))=\zeta\mu(m)$;
 the map $\zeta$ is multiplicative, since $\zeta(a)\zeta(b)=\partial(0)=0=\zeta(ab)$ for $a,b\in\ZZ_{\MM}(\A)$; and the actions are preserved,
  since ${}^{\zeta(a)}m={}^am=0$ and $m^{\zeta(a)}=m^a=0$ for $a\in\ZZ_{\MM}(\A)$ and $m\in\Ann_{\MM}(\A)$. Clearly, $\inc$ and $\zeta$ are injective.

Moreover, the following diagram is commutative
\[	
\begin{tikzcd}[column sep=0.3cm]
	0 \arrow[r] &  \pi_1 \big(\ZZ^{\text{Huq}}(\MM,\A,\mu) \big) \arrow[d,equal]   \arrow[r] & \Ann_{\MM}(\A) \arrow[d,hook,"\inc"] \arrow[r, "\mu|_{\Ann_{\MM}(\A)}"]  
	& \ZZ_{\MM}(\A)   \arrow[d, "\zeta"]   \arrow[r] &  \pi_0 \big(\ZZ^{\text{Huq}}(\MM,\A,\mu) \big) \arrow[d,  "\overline{\zeta}"]   \arrow[r] & 0   \\
	0 \arrow[r] &  \pi_1 \big(\ZC(\MM,\A,\mu)\big)  \arrow[r]   & \MM  \arrow[r, "\partial"] & \T(\MM,\A,\mu)  \arrow[r]   &  \pi_0 \big(\ZC(\MM,\A,\mu)\big)   \arrow[r] & 0 
\end{tikzcd} 
\]
Additionally, $\pi_1 \big(\ZZ^{\text{Huq}}(\MM,\A,\mu) \big)= \Ann_{\MM}(\A) \cap \Ker \mu=\pi_1 \big(\ZC(\MM,\A,\mu)\big)$  and 
$\overline{\zeta}$ is injective because $\zeta\big(\ZZ_{\MM}(\A)\big) \cap \Im \partial=\{(0,0, \mu(m)) \mid \mu(m) \in  \ZZ_{\MM}(\A) \ \text{and}  \ m \in \Ann_{\MM}(\A)\} \cong \mu\big(\Ann_{\MM}(\A)\big)$.

\end{proof}
\begin{remark}
	In general, the morphism $(\inc,\zeta)$ is not a weak equivalence, since $\overline{\zeta}$ need not be an isomorphism.
For instance, let $\A$ be the algebra of Example~\ref{E:CM}~(3) and let $\MM=\K f$ be its ideal, regarded as an abelian associative algebra and as an $\A$-bimodule via the product of $\A$.
 Since ${}^ef=ef=f$ and $\ZZ(\A)=0$, we get $\Ann_{\MM}(\A)=0=\ZZ_{\MM}(\A)$, so $\ZZ^{\text{Huq}}(\MM,\A,0)=(0,0,0)$. On the other hand,
  the right action of $\A$ on $\MM$ is trivial ($fe=0=f^2$), so a left multiplication $\lam\colon\A\to\K f$ satisfies $\lam(e)=\lam(e^2)=\lam(e)^e=0$
   and $\lam(f)=\lam(ef)=\lam(e)^f=0$, whereas any linear map $\rho\colon\A\to\K f$ is a right multiplication, and \eqref{Mbimul} holds trivially. Thus
	\[
	\T(\MM,\A,0)=\Mul(\A,\MM)=\{(0,\rho,0)\mid \rho\colon \A\to\K f \ \text{linear}\}\cong\K^2,
	\]
	while $\partial(f)=(0,\RR_f,0)$ with $\RR_f(e)=f$ and $\RR_f(f)=0$. Therefore $\pi_0\big(\ZC(\MM,\A,0)\big)\cong\K$, generated by the class of $(0,\rho,0)$
	 with $\rho(e)=0$ and $\rho(f)=f$, and $\overline{\zeta}\colon 0\to\K$ is not an isomorphism.
\end{remark}


\section{A Drinfeld-type interpretation of the center} \label{S:Drin}
It is well known that the category of crossed modules of associative algebras is equivalent to the category of internal categories in associative algebras \cite[Theorem~3.5]{Emz}
 (see also \cite[Proposition~1.1]{FFL1}), or shortly, the category of associative 2-algebras. Given $(\MM,\A,\mu)$ a crossed module of associative algebras,
  we denote by $\Cat(\MM,\A,\mu)$ its associated internal category whose objects are the elements of $\A$, 
  and morphisms $a\xrightarrow{m \,}b$ are labeled by the elements $m\in \MM$ satisfying $\mu(m)+a=b$. The identity arrows are labeled by $0\in\MM$, whereas the composition is given by addition of labels. 

In this section, we provide a categorical interpretation of the center constructed in Definition~\ref{D:center}, inspired by the notion of the Drinfeld center \cite{Ma} (see also \cite[Section~7.13]{EGNO}). 

Since $\Cat(\MM,\A,\mu)$ is an internal category, the multiplication in $\A$ induces a bifunctor
 \[\cdot\colon \Cat(\MM,\A,\mu)\times \Cat(\MM,\A,\mu)\longrightarrow\Cat(\MM,\A,\mu)\]
  defined on objects by $(a,b)\mapsto ab$, and on morphisms by
\[\big(a\xrightarrow{m\,}b\big)\cdot\big(a'\xrightarrow{m' \,}b'\big)\coloneqq \big(aa'\xrightarrow{mm'+{}^{a}m'+m^{a'}}bb'\big).\]
In particular, for any object $c$, we have $\id_c \cdot \big(a\xrightarrow{m\,}b\big)=\big(ca\xrightarrow{{}^{c}m}cb\big)$ and $\big(a\xrightarrow{m\,}b\big) \cdot\id_c =\big(ac\xrightarrow{m^c}bc\big)$.

The label of the product of two morphisms is the first component of the product in the semidirect product algebra $\MM\rtimes\A$ of Section~\ref{S:Prel}.
 The associativity of this product ensures that the bifunctor $\cdot$ is strictly associative. However, since associative algebras considered in this article are not necessarily unital,
  $(\Cat(\MM,\A,\mu),\cdot)$ is in general not a strict monoidal category, only \emph{semimonoidal} (see~\cite{JK}).

Let us denote by $\mathbf{0}$ the zero endofunctor of $\Cat(\MM,\A,\mu)$, which assigns to any object $a\in\A$ the zero element $0\in\A$, and consider the endofunctors $a\cdot (-)$ and $(-)\cdot a$.

Given $(\lambda,\rho,a)\in \T(\MM,\A,\mu)$, for any $b\in \A$ we define the following two families of morphisms in $\Cat(\MM,\A,\mu)$:
\[\boldsymbol{\lambda}_{b}\colon 0\xrightarrow{\lambda(b)}ab,\qquad \boldsymbol{\rho}_{b}\colon 0\xrightarrow{\rho(b)}ba,\]
which are well-defined thanks to conditions \eqref{Z1} and \eqref{Z2}. Furthermore, axioms~\eqref{Z3} and \eqref{Z4} imply that the above families determine natural transformations 
\[\boldsymbol{\lambda} \colon \mathbf{0} \Longrightarrow a\cdot(-),\qquad \boldsymbol{\rho} \colon\mathbf{0} \Longrightarrow (-)\cdot a,\]
respectively. Moreover, in terms of the preceding natural transformations, 
\begin{align}\label{C1}\tag{C1}
	\lambda\in\LM(\A,\MM)&\iff \boldsymbol{\lambda}_{bc}=\boldsymbol{\lambda}_{b}\cdot \id_{c},\\
	\label{C2}\tag{C2}\rho\in\RM(\A,\MM)&\iff \boldsymbol{\rho}_{bc}= \id_{b}\cdot \boldsymbol{\rho}_{c},\\
	\label{C3}\tag{C3}\eqref{Mbimul}&\iff \id_{b}\cdot \boldsymbol{\lambda}_{c}=\boldsymbol{\rho}_{b}\cdot \id_{c}.
\end{align}

The previous remarks motivate the following definition:
\begin{definition}
	Let $(\MM,\A,\mu)$ be an associative crossed module. The \emph{center of an associative 2-algebra} $\ZC\big(\Cat(\MM,\A,\mu)\big)$ 
	is the category whose objects are triples $(\boldsymbol{\lambda} ,\boldsymbol{\rho},a)$, where $a\in\A$ and
	\[\boldsymbol{\lambda}  \colon \mathbf{0} \Longrightarrow a\cdot(-),\qquad \boldsymbol{\rho} \colon\mathbf{0} \Longrightarrow (-)\cdot a\]
	are natural transformations with components determined by linear maps $\boldsymbol{\lambda} _{b}\colon 0\xrightarrow{\lambda(b)}ab$ and $\boldsymbol{\rho}_{b}\colon 0\xrightarrow{\rho(b)}ba$,
	 satisfying \eqref{C1}, \eqref{C2} and \eqref{C3}.
	
	A morphism $(\boldsymbol{\lambda} ,\boldsymbol{\rho},a)\rightarrow (\boldsymbol{\lambda} ',\boldsymbol{\rho}',b)$ in $\ZC(\Cat(\MM,\A,\mu))$ is a morphism
	 $a\xrightarrow{m}b$ in $\Cat(\MM,\A,\mu)$ such that, for any $c\in\A$, the diagrams
	
	\noindent\begin{minipage}{0.48\linewidth}
		\[\xymatrix{&0\ar[rr]^-{\boldsymbol{\lambda} _{c}}\ar[rrd]_-{\boldsymbol{\lambda} '_{c}} & &ac\ar[d]^-{m\cdot \id_{c}}\\& & &bc}\]
	\end{minipage}\hfill
	\begin{minipage}{0.48\linewidth}
		\begin{equation}\label{morcat}
			\begin{aligned}
				\xymatrix{&0\ar[rr]^-{\boldsymbol{\rho}_{c}}\ar[rrd]_-{\boldsymbol{\rho}'_{c}} & &ca\ar[d]^-{\id_{c}\cdot m}\\& & &cb}
			\end{aligned}
		\end{equation}
	\end{minipage}
	commute.
\end{definition}
 
The following theorem identifies this category with the internal category associated with the center crossed module constructed in Definition~\ref{D:center}.

\begin{theorem}
	Let $(\MM,\A,\mu)$ be an associative crossed module. There exists a canonical isomorphism of categories between $\ZC(\Cat(\MM,\A,\mu))$ and $\Cat(\ZC(\MM,\A,\mu))=\Cat(\MM,\T(\MM,\A,\mu),\partial)$.
\end{theorem}
\begin{proof}
	Let $(\boldsymbol{\lambda} ,\boldsymbol{\rho},a)\in \ZC\big(\Cat(\MM,\A,\mu)\big)$. Writing $\lambda(b)$ and 
	$\rho(b)$ for the labels of $\boldsymbol{\lambda} _{b}$ and $\boldsymbol{\rho}_{b}$, respectively,
	 we have that $\mu\lambda(b)=ab$ and $\mu\rho(b)=ba$. Moreover, the naturality of $\boldsymbol{\lambda} $
	  and $\boldsymbol{\rho}$ yields $\lambda\mu(m)={}^{a}m$ and $\rho\mu(m)=m^{a}$.
	   Therefore, taking into account the equivalences given in \eqref{C1}--\eqref{C3},
	    we obtain that $(\lambda,\rho)\in\Mul(\A,\MM)$, and then $(\lambda,\rho, a)\in \T(\MM,\A,\mu)$. 
	    Conversely, every element of $\T(\MM,\A,\mu)$ determines a unique object in $\ZC\big(\Cat(\MM,\A,\mu)\big)$ in this way.
	
	Now, let $(\boldsymbol{\lambda} ,\boldsymbol{\rho},a)\xrightarrow{m \,}(\boldsymbol{\lambda} ',\boldsymbol{\rho}',b)\in \ZC\big(\Cat(\MM,\A,\mu)\big)$.
	 Since $a\xrightarrow{m\, }b\in \Cat(\MM,\A,\mu)$, we have that $b=a+\mu(m)$. Moreover, compatibilities given in diagrams \eqref{morcat}
	are equivalent to
	\[\lambda'=\lambda+\LL_{m},\qquad \rho'=\rho+\RR_{m}.\]
	Consequently, $(\lambda',\rho',b)=(\lambda,\rho,a)+\partial(m)$, which is exactly the condition defining a morphism in $\Cat(\ZC(\MM,\A,\mu))$.
	 Since the correspondence preserves identities and composition, it defines the required canonical isomorphism of categories.
\end{proof}

Once we have shown that $\ZC\big(\Cat(\MM,\A,\mu)\big)\cong \Cat\big(\ZC(\MM,\A,\mu)\big)$, a natural question that arises is whether the braiding on 
$\ZC(\MM,\A,\mu)$ defined in \eqref{braiding_center} induces a categorical braiding on $\Cat\big(\ZC(\MM,\A,\mu)\big)$. 
Let $(\lambda,\rho,a),(\lambda',\rho',b)\in\T(\MM,\A,\mu)$ and consider the morphisms in $\Cat(\MM,\A,\mu)$ given by $0\xrightarrow{\lambda(b)}ab$ and $0\xrightarrow{\lambda'(a)}ba$. Since
\[(\lambda,\rho,a)(\lambda',\rho',b)=\partial\lambda(b), \qquad(\lambda',\rho',b)(\lambda,\rho,a)=\partial\lambda'(a),\]
the linearity of $\partial$ yields \[\partial\{(\lambda,\rho,a),(\lambda',\rho',b)\}=\partial(\lambda(b)-\lambda'(a))=(\lambda,\rho,a)(\lambda',\rho',b)-(\lambda',\rho',b)(\lambda,\rho,a).\]
As a consequence, $\{(\lambda,\rho,a),(\lambda',\rho',b)\}$ determines a morphism
\[\beta_{(\lambda,\rho,a),(\lambda',\rho',b)}\colon (\lambda',\rho',b)(\lambda,\rho,a)\xrightarrow{\{(\lambda,\rho,a),(\lambda',\rho',b)\}}(\lambda,\rho,a)(\lambda',\rho',b)\]
in $\Cat\big(\ZC(\MM,\A,\mu)\big)$. Note also that every morphism of an associative 2-algebra is invertible. Hence, $\beta_{(\lambda,\rho,a),(\lambda',\rho',b)}$ is an isomorphism.

In the next result we will show that $\beta$ defines a braiding for $\Cat\big(\ZC(\MM,\A,\mu)\big)$.
\begin{proposition}
	The family $\beta$ constructed above defines a natural transformation 
	\[\beta\colon \cdot^{\text{\emph{op}}} \Longrightarrow \cdot\] satisfying the strict hexagon identities:
	\begin{gather}
	\label{braid1}	\beta_{(\lambda,\rho,a)(\lambda',\rho',b),(\lambda'',\rho'',c)}=(\id_{(\lambda,\rho,a)}\cdot \beta_{(\lambda',\rho',b),(\lambda'',\rho'',c)})\circ (\beta_{(\lambda,\rho,a),(\lambda'',\rho'',c)}\cdot \id_{(\lambda',\rho',b)}),\\\label{braid2}
		\beta_{(\lambda,\rho,a),(\lambda',\rho',b)(\lambda'',\rho'',c)}=(\beta_{(\lambda,\rho,a),(\lambda',\rho',b)}\cdot \id_{(\lambda'',\rho'',c)})\circ ( \id_{(\lambda',\rho',b)}\cdot \beta_{(\lambda,\rho,a),(\lambda'',\rho'',c)}).
	\end{gather}
	
	In other words, $\big(\Cat(\ZC(\MM,\A,\mu)),\cdot,\beta\big)$ is a braided semimonoidal category.
\end{proposition}
\begin{proof}
	Let us begin by showing the naturality of $\beta$ in the first variable. The same argument applies for the second one, so these details will be omitted.
	 Let $(\lambda,\rho,a)\xrightarrow{m\,}(\lambda+\LL_{m},\rho+\RR_{m},a+\mu(m))$ be a morphism in $\Cat\big(\ZC(\MM,\A,\mu)\big)$ and fix
	  $(\lambda',\rho',b)\in \Cat\big(\ZC(\MM,\A,\mu)\big)$. We have to show that 
	\[\beta_{(\lambda+\LL_{m},\rho+\RR_{m},a+\mu(m)),(\lambda',\rho',b)}\circ (\id_{(\lambda',\rho',b)}\cdot m)=(m\cdot  \id_{(\lambda',\rho',b)})\circ \beta_{(\lambda,\rho,a),(\lambda',\rho',b)}\]
	which is equivalent to 
	\[\{(\lambda+\LL_{m},\rho+\RR_{m},a+\mu(m)),(\lambda',\rho',b)\}+{}^{b}m=m^{b}+\{(\lambda,\rho,a),(\lambda',\rho',b)\}.\]
	Indeed, the previous identity follows by:
	\begin{align*}
		\{(\lambda+\LL_{m},\rho+\RR_{m},a+\mu(m)),(\lambda',\rho',b)\}&=(\lambda+\LL_{m})(b)-\lambda'(a+\mu(m))\\
		&=\{(\lambda,\rho,a),(\lambda',\rho',b)\}+m^{b}-\lambda'\mu(m)\\
		&\stackrel{\eqref{Z3}}{=}\{(\lambda,\rho,a),(\lambda',\rho',b)\}+m^{b}-{}^{b}m.
	\end{align*}
	
	In what follows, we will show the hexagon identity \eqref{braid1} to hold, whereas \eqref{braid2} is analogous. Note that, by using \eqref{eq:prod}, \eqref{braid1} can be rewritten as
	\[\beta_{\partial\lambda(b),(\lambda'',\rho'',c)}=(\id_{(\lambda,\rho,a)}\cdot \beta_{(\lambda',\rho',b),(\lambda'',\rho'',c)})\circ (\beta_{(\lambda,\rho,a),(\lambda'',\rho'',c)}\cdot \id_{(\lambda',\rho',b)}),\]
	which is equivalent to
	\[\{\partial\lambda(b),(\lambda'',\rho'',c)\}={}^{a}\{(\lambda',\rho',b),(\lambda'',\rho'',c)\}+\{(\lambda,\rho,a),(\lambda'',\rho'',c)\}^{b}.\]
	The equality above results from
	\begin{align*}
	& \quad \	{}^{a}\{(\lambda',\rho',b),(\lambda'',\rho'',c)\}+\{(\lambda,\rho,a),(\lambda'',\rho'',c)\}^{b}\\
		&={}^{a}(\lambda'(c)-\lambda''(b))+(\lambda(c)-\lambda''(a))^{b}={}^{a}\lambda'(c)-{}^{a}\lambda''(b)+\lambda(c)^{b}-\lambda''(a)^{b}\\
		&\stackrel{\eqref{Mlemu},\eqref{Z3}}{=}\lambda\mu\lambda'(c)-\lambda\mu\lambda''(b)+\lambda(cb)-\lambda''(ab)\stackrel{\eqref{Z1}}{=}\lambda(bc)-\lambda''(ab)
	\end{align*}
	and, on the other hand, 
	\begin{align*}
		\{\partial\lambda(b),(\lambda'',\rho'',c)\}&=\{(\LL_{\lambda(b)},\RR_{\lambda(b)},ab),(\lambda'',\rho'',c)\}\\
		&=\LL_{\lambda(b)}(c)-\lambda''(ab)=\lambda(b)^{c}-\lambda''(ab)\stackrel{\eqref{Mlemu}}{=}\lambda(bc)-\lambda''(ab).\qedhere
	\end{align*}
\end{proof}

The preceding construction provides a Drinfeld-type interpretation of the center $\ZC(\MM,\A,\mu)$. In fact, each $(\lambda,\rho,a)\in\T(\MM,\A,\mu)$ determines a natural transformation 
\[\alpha\colon (-)\cdot a \Longrightarrow a\cdot (-),\]
given by the family of morphisms $\alpha_{b}\colon ba\xrightarrow{\lambda(b)-\rho(b)}ab$, $b\in \A$. However, the half-braiding $\alpha=\boldsymbol{\lambda}-\boldsymbol{\rho}$ does not necessarily determine the pair $(\lambda,\rho)$ uniquely. Thus, our construction retains more information than the corresponding natural isomorphism $\alpha$ and should be regarded as a Drinfeld-type analogue rather than an identification with the ordinary Drinfeld center.


\section{Homotopical interpretation of the center} \label{S:homoto}

Regard an associative crossed module $(\MM,\A,\mu)$ as a $2$-term chain complex $\MM\xrightarrow{\mu}\A$
concentrated in degrees $1$ and $0$.

A chain map  $f=(f_1,f_0)\colon  (\MM,\A,\mu) \to  (\MM',\A',\mu')$  is a pair of linear maps such that $f_0\mu=\mu' f_1$. A \emph{chain homotopy}  between two chain maps $f=(f_1,f_0)$ and  $g=(g_1,g_0)$,
denoted by $f \simeq g$, is a
linear map $h \colon \A \to \MM'$ such that $f_1 - g_1= h \mu$ and   $f_0 - g_0=  \mu' h$. 
A chain endomorphism   $f=(f_1,f_0)$ is called  \emph{null-homotopic} if $f \simeq 0$, that is,  if there exists a linear map $h \colon \A \to \MM$ such that  $f_1 = h \mu$ and   $f_0=  \mu h$. In this  case $h$ is called a  \emph{null-homotopy} of $f$.

\begin{proposition}\label{prop:homotopy}
	Let $(\MM,\A,\mu)$ be an associative crossed module and $a\in\A$. Consider the two pairs of linear maps (by a slight abuse of notation, we keep the names $\LL_a$ and $\RR_a$)
	\[
	\LL_a=\bigl(\MM \xrightarrow{m\mapsto\lact{a}{m}} \MM,\ \A \xrightarrow{ \, a\cdot \ } \A\bigr),
	 \qquad \RR_a=\bigl(\MM \xrightarrow{m\mapsto\ract{m}{a}} \MM,\ \A \xrightarrow{ \, \cdot a  \ } \A \bigr).
	\]
	Then the following assertions hold:
	\begin{enumerate}
		\item $\LL_a$ and $\RR_a$  are chain endomorphisms of $(\MM,\A,\mu)$.
		\item A linear map $\lambda \colon\A\to\MM$ is a null-homotopy of $\LL_a$ if and
		only if  \eqref{Z1}  and \eqref{Z3}  hold. Similarly,  a linear map $\rho \colon\A\to\MM$ is a null-homotopy of $\RR_a$ if and only if  \eqref{Z2}
		and \eqref{Z4} hold. 
		\item Consequently
		\[
		\T(\MM,\A,\mu)=\Bigl\{(\lambda,\rho,a)\ \Bigm|\
		\begin{array}{l}
			\lambda \ \text{is a null-homotopy of}\ \LL_a,\\
			\rho\ \text{is a null-homotopy of} \ \RR_a, \ \text{and}\\
			\text{the pair}\ (\lambda,\rho) \ \text{is a bimultiplication}
		\end{array}\Bigr\}.
		\]
	\end{enumerate}
\end{proposition}
\begin{proof}
	The proof follows directly from the definitions.
\end{proof}

\begin{remark}
The center of an associative algebra in the categorical sense is its annihilator, i.e. it is characterized  by
$\LL_a=0=\RR_a$. In contrast, the center of an associative  crossed module can be described as the crossed module of elements of $\A$ equipped with coherent null-homotopies.
\end{remark}



\begin{thebibliography}{99}
	
\bibitem{AH} U.~E. Arslan  and S. H\"urmetli, Bimultiplications and annihilators of crossed modules in
	associative algebras,  \emph{Journal of New Theory} \textbf{35} (2021),  72--90.
	
\bibitem{BG} D. Bourn and M. Gran, Central extensions in semi-abelian categories, \emph{J. Pure Appl. Algebra} \textbf{175} (1-3) (2002),  31--44.
	
\bibitem{CC}
P. Carrasco and A.~M. Cegarra, A Dold--Kan theorem for simplicial Lie algebras, \emph{Theory Appl. Categ.} \textbf{32} (2017), 1165--1212.

\bibitem{CDL}
J.~M. Casas, T.  Datuashvili and M. Ladra, 
Universal strict general actors and actors in categories of interest,
\emph{Appl. Categ. Structures} \textbf{18} (1) (2010),  85--114.

\bibitem{CaLa}
J.~M. Casas and M. Ladra, The actor of a crossed module in Lie algebras, \emph{Comm. Algebra} \textbf{26} (7) (1998), 2065--2089.

\bibitem{DL} P.~Dedecker and A. S.-T. Lue, 
A nonabelian two-dimensional cohomology for associative algebras,
\emph{Bull. Amer. Math. Soc.} \textbf{72} (1966), 1044--1050.

\bibitem{El}
G.~J. Ellis,  \emph{Crossed modules and their higher dimensional analogues}, PhD thesis, University of Wales (Bangor), 1984.

 \bibitem{EGNO} P.~I. Etingof, S. Gelaki, D. Nikshych and V. Ostrik, \emph{Tensor categories}, Mathematical Surveys and Monographs \textbf{205}, Amer. Math. Soc., Providence, RI, 2015.

 \bibitem{Ev} C.~J. Everett, An extension theory for rings, \emph{Amer. J. Math.} \textbf{64} (1942), 363--370.

 \bibitem{FFL1}
  A.~Fern\'andez-Fari\~na and M.~Ladra,   Braiding for categorical algebras and crossed modules of algebras I: Associative and Lie algebras, \emph{J. Algebra Appl.} \textbf{19} (9) (2020), 2050176, 30 pp.

  \bibitem{Ger} M.~Gerstenhaber, On the deformation of rings and algebras. II,
\emph{Ann. of Math. (2)} \textbf{84} (1966), 1--19.

 \bibitem{Hoc} G.~Hochschild, Cohomology and representations of associative algebras, \emph{Duke Math. J.} \textbf{14} (1947), 921--948.

 \bibitem{Huq} S. A. Huq,  Commutator, nilpotency, and solvability in categories, \emph{Quart. J. Math.
  Oxford Ser. (2)} \textbf{19}  (1968), 363--389.
 
\bibitem{J} G. Janelidze, Internal crossed modules, \emph{Georgian Math. J.} \textbf{10} (1) (2003), 99--114.

  \bibitem{JMT} G. Janelidze, L. M\'arki and W. Tholen, Semi-abelian categories, \emph{J. Pure Appl. Algebra} \textbf{168} (2-3) (2002),  367--386.
  
  \bibitem{JK}
  A. Joyal and J. Kock, \emph{Weak units and homotopy 3-types}, Categories in algebra, geometry and mathematical physics, Contemporary Mathematics \textbf{431}  (2007),  257--276.
  
  \bibitem{JS}
  A. Joyal and R.~H. Street, Tortile Yang-Baxter operators in tensor categories, \emph{J. Pure Appl. Algebra} \textbf{71} (1) (1991),  43--51.
  
    \bibitem{KL} 
  C. Kassel and J.-L. Loday, Extensions centrales d'alg\`ebres de Lie, \emph{Ann. Inst. Fourier (Grenoble)} \textbf{32} (4) (1982),  119--142. 

\bibitem{Emz} E. Khmaladze, On associative and Lie 2-algebras, \emph{Proc. A. Razmadze Math. Inst.} \textbf{159} (2012), 57--64.

\bibitem{LL}
R.~Lavendhomme  and Th. Lucas,  On modules and crossed modules, \emph{J. Algebra} \textbf{179} (3) (1996),  936--963.

 \bibitem{Mac} S. Mac~Lane, Extensions and obstructions for rings, \emph{Illinois J. Math.} \textbf{2} (1958), 316--345.
    
  \bibitem{Ma}
   S. Majid, Representations, duals and quantum doubles of monoidal categories, \emph{Rend. Circ. Mat. Palermo (2)} Suppl.  \textbf{26} (1991), 197--206.
  
 \bibitem{Nor}
K. Norrie, Actions and automorphisms of crossed modules, \emph{Bull. Soc. Math. France} \textbf{118} (2) (1990),  129--146.

 
   
 \bibitem{Pim1}
 M.~Pirashvili, On the centre of crossed modules of groups, \emph{J. Math. Sci. (N.Y.)} \textbf{290} (5) (2025),  587--597. 

 \bibitem{Pim2}
  M.~Pirashvili, On the centre of crossed modules of Lie algebras, \emph{J. Lie Theory} \textbf{35} (2) (2025),   359--375.


  
 \bibitem{UCA}
  E.~\"O. Uslu, S. \c Cetin and A.~F. Arslan, On crossed modules in modified categories of interest, \emph{Math. Commun.} \textbf{22} (1) (2017), 103--119.
  
   \bibitem{Wa}
  F.~Wagemann, \emph{Crossed modules}, De Gruyter Studies in Mathematics, 82, De Gruyter, Berlin, 2021. 
  
   \bibitem{Wh}
  J. H. C.~Whitehead, Combinatorial homotopy II, \emph{Bull. Amer. Math. Soc.} \textbf{55} (1949), 453--496.

 \end{thebibliography}
\end{document}